\documentstyle[12pt,amssymb,amscd]{amsart}  

\theoremstyle{remark}

\theoremstyle{definition}

\numberwithin{equation}{subsection}

\newcommand{\bbA}{{\Bbb A}}

\newcommand{\bbC}{{\Bbb C}}

\newcommand{\bbH}{{\Bbb H}}

\newcommand{\CC}{{\Bbb C}}

\newcommand{\bbR}{{\Bbb R}}

\newcommand{\bbZ}{{\Bbb Z}}
\newcommand{\cF}{{\cal F}}
\newcommand{\cG}{{\cal G}}
\newcommand{\cO}{{\cal O}}

\newcommand{\cL}{{\cal L}}
\newcommand{\cM}{{\cal M}}

\newcommand{\cD}{{\cal D}}

\newcommand{\cA}{{\cal A}}

\newcommand{\cE}{{\cal E}}

\newcommand{\cU}{{\cal U}}

\newcommand{\cT}{{\cal T}}

\newcommand{\cH}{{\cal H}}

\newcommand{\Lift}{\operatorname{Lift}}
\newcommand{\Coker}{\operatorname{Coker}}
\newcommand{\Aut}{\operatorname{Aut}}
\newcommand{\DR}{\operatorname{DR}}

\newcommand{\Ker}{\operatorname{Ker}}

\newcommand{\Hom}{\operatorname{Hom}}
\newcommand{\End}{\operatorname{End}}

\newcommand{\id}{\operatorname{id}}
\newcommand{\Tr}{\operatorname{Tr}}
\newcommand{\shHom}{\underline{\operatorname{Hom}}}
\newcommand{\shEnd}{\underline{\operatorname{End}}}

\newcommand{\Smbl}{\operatorname{Smbl}}
\newcommand{\Ext}{\operatorname{Ext}}

\newcommand{\ch}{\operatorname{ch}}
\newcommand{\rk}{\operatorname{rk}}

\newcommand{\ad}{\operatorname{ad}}

\newcommand{\Lie}{\operatorname{Lie}}

\newcommand{\Gg}{{\frak g}}

\newcommand{\Ad}{\operatorname{Ad}}

\newcommand{\Td}{\operatorname{Td}}
\newcommand{\Det}{\operatorname{Det}}
\newcommand{\Hoch}{\operatorname{Hoch}}
\newcommand{\gr}{\operatorname{gr}}
\newcommand{\GA}{{\mathfrak A}}
\newcommand{\Ob} {\operatorname{Ob}}
\newcommand{\Mor}{\operatorname{Mor}}
\newcommand{\Res}{\operatorname{Res}}

\begin{document}

\title{Riemann-Roch for real varieties}

\author{Paul Bressler}
\email{bressler@@math.arizona.edu}
\address{University of Arizona}
\author{Mikhail Kapranov}
\email{mikhail.kapranov@@yale.edu}
\address{Yale University}
\author{Boris Tsygan}
\email{tsygan@@math.northwestern.edu}
\address{Northwestern University}
\author{Eric Vasserot}
\email{vasserot@@math.jussieu.fr}
\address{Universit\'{e} Paris-7}
\maketitle

 \begin{center}
\vspace*{.7cm}
{\it To Yuri Ivanovich Manin on his 70th birthday}
\end{center}

\vskip 1cm

\centerline {\bf Introduction}

\vskip 1cm

\noindent {\bf (0.1)} Let $\Sigma$ be an oriented
 real analytic manifold of dimension $d$ and $X$ be a complex
envelope of $\Sigma$, i.e., a complex manifold of the same dimension containing
$\Sigma$ as a totally real submanifold. Then, (real) geometric objects
on $\Sigma$ can be viewed as (complex) geometric objects on $X$ involving
cohomology classes of degree $d$. For example, a $C^\infty$-function $f$
on $\Sigma$ can be considered as a section of ${\mathcal B}_\Sigma$, the
sheaf of hyperfunctions on $\Sigma$ which, according to Sato, can be
defined as
$${\mathcal B}_\Sigma = \underline{H}^d_\Sigma(\cO_X), \leqno (0.1.1)
$$
where $\cO_X$ is the sheaf of holomorphic functions. 
So $f$ can be viewed as a class in $d$th local cohomology. 

More generally, the equality (0.1.1) suggests that various results of
holomorphic geometry on $X$ should have consequences for the purely real geometry on $\Sigma$,
consequences that
involve raising the cohomological degree by $d$. The goal of this paper is to
investigate the consequences of one such result, the Grothendieck-Riemann-Roch theorem
(GRR). 

\vskip .3cm

\noindent {\bf (0.2)} Let $p: X\to B$ be a smooth proper morphism
of complex algebraic manifolds. We denote the fibers of $p$ by $X_b=p^{-1}(b)$ and assume
them to be of dimension $d$. If $\mathcal{E}$ is an algebraic vector bundle on $X$, the GRR
theorem says that
$$ch_m(Rp_*({\mathcal{E}})) = \int_{X/B} \biggl[ ch({\mathcal E})\cdot \Td({\mathcal T}_{X/B})\biggr]_{2m+2d}
\quad \in \quad  H^{2m}(B,
\mathbb{C}). \leqno (0.2.1)$$
Here $\int_{X/B}: H^{2m+2d}(X, {\mathbb C})\to H^{2m}(B, {\mathbb C})$ is the cohomological
direct image (integration over the fibers of $p$). 

In the case $m=1$ the class on the left comes from the class, in the Picard group of $B$,
of the determinantal line bundle $\det (Rp_*{\mathcal E})$ whose fiber, at a generic
point $b\in B$, is
$$\det\, H^\bullet(X_b, {\mathcal E}) \quad =\quad \bigotimes_i \left( \Lambda^{max} \, H^i(X_b, {\mathcal E})
\right) ^{\otimes \, (-1)^i}. \leqno (0.2.2)$$
 Deligne [D1] posed the problem of describing $\det(Rp_*{\mathcal E})$
in a functorial way  as  a refinement of  GRR  for $m=1$.
This problem makes sense already for the case $B=pt$ when we have to describe the
1-dimensional vector space (0.2.2) as a functor of $\mathcal E$. 
Deligne  solved this problem for a family of curves and further results have been obtained in [E]. 

\vskip .3cm

\noindent {\bf (0.3)} To understand the real counterpart of (0.2.1), assume first that
$B=pt$, so $X=X_{pt}$ and let $\Sigma\subset X$ be as in (0.1). Denote by $E$ the restriction
of $\mathcal{E}$ to $\Sigma$ and by $C^\infty_\Sigma(E)$ the sheaf of its $C^\infty$ sections.
Then, similarly to (0.1.1) we have the embedding
$$C^\infty_\Sigma(E)\subset \underline{H}^d_\Sigma({\mathcal E}).
\leqno (0.3.1)$$
Assume further that $d=1$, so $X$ is an algebraic curve, and that $\Sigma$ is a small
circle in $X$ cutting it into two pieces: $X_+$ (a small disk) and $X_-$. Let $\mathcal{E}_\pm = 
{\mathcal E}|_{X_\pm}$. We are then in the situation of the Krichever correspondence [PS].
Namely, the space $\Gamma(E)$ of $L^2$-sections has a canonical polarization in the sense
of Pressley and Segal [PS] and therefore possesses a determinantal gerbe $\Det \, \Gamma(E)$.
The latter is a category with every Hom-set  made into a $\mathbb{C}^*$-torsor (a 1-dimensional
vector space with zero deleted). The extensions ${\mathcal E}_\pm$ of $E$ to $X_\pm$ define
two objects $[\mathcal {E}_\pm]$ of this gerbe, and
$$\det\, H^\bullet(X,{ \mathcal E})  = \Hom_{\Det\, \Gamma(E)} ([{\mathcal E}_+], [{\mathcal E}_-]).
\leqno (0.3.2)$$
The real counterpart of the problem of describing the $\mathbb{C}^*$-torsor $\det\, H^\bullet(X,
{\mathcal E})$ is then the problem of describing the gerbe $\Det\, \Gamma(E)$. If we now have a family 
$p: X\to B$ as before (with $d=1$), equipped with a subfamily of circles $q: \Sigma\to B$,\
$\Sigma\subset X$, then we have an $\mathcal{O}^*_B$-gerbe $\Det \, q_*(E)$ which, according to the
the classification of gerbes [Bre],  has a class in $H^2(B, \mathcal{O}^*_B)$. The latter group
maps naturally to $H^3(B, \mathbb{Z})$ and in fact can be identified with the Deligne
cohomology group $H^3(B, {\mathbb Z}_{ D}(1))$, see [Bry]. The Real Riemann-Roch for a
circle fibration describes the above class (modulo 2-torsion) as
$$[\Det\, q_*(E)] = \int_{\Sigma/B} ch_2(E)\quad  \in\quad H^3(B, {\mathbb Z}_{\mathcal D}(2))\otimes
{\mathbb Z}\left[ {1\over 2}\right].
\leqno (0.3.3)$$
Here $\int_{\Sigma/B}: H^4(\Sigma, {\mathbb Z}_{ D}(2))\to H^3(B, {\mathbb Z}_{ D}(1))$
is the direct image in Deligne cohomology.  Note the absense of the characteristic classes
of ${\mathcal T}_{\Sigma/B}$ (they are 2-torsion for a real rank one bundle). If one is
interested in the image of the determinantal class in $H^3(B, \mathbb{Z})$, then one can 
understand the RHS of the above formula in the purely topological sense.

\vskip .2cm

Both sides of (0.3.3) do not involve anything other than $q:\Sigma\to B$ and a vector
bundle $E$ on $\Sigma$ (equipped with CR-structures coming from the embeddings into
$X,\mathcal{E}$). One has a similar result for any $C^\infty$ circle fibration
(no CR structure) and any $C^\infty$ complex bundle $E$ on $\Sigma$. In this case we 
get a gerbe with lien $C_B^{\infty *}$, the sheaf of invertible complex valued $C^\infty$-functions
on $B$ and its class lies in $H^2(B, C_B^{\infty *}) = H^3(B, \mathbb{Z})$. 
It is this, purely $C^\infty$ setting, that we adopt and generalize in the
present paper. 

\vskip .3cm

\noindent {\bf (0.4)} Let now $\Sigma$ be a compact oriented $C^\infty$-manifold of arbitrary
dimension $d$ and $E$ a $C^\infty$ complex vector bundle on $\Sigma$. One expects that the space
$\Gamma(E)$ should have some kind of $d$-fold polarization, giving rise to a
``determinantal $d$-gerbe'', $\Det\, \Gamma(E)$. This structure is rather clear when
$\Sigma$ is a 2-torus but in general the theory of higher gerbes is not fully
developed.	In any case one expects that a $C^\infty$ family of such gerbes over a base $B$
gives a class in $H^{d+1}(B, C_B^{\infty *}) = H^{d+2}(B, \mathbb{Z})$. In this paper we
consider a $C^\infty$ family $q:\Sigma\to B$ of relative dimension $d$ and a $C^\infty$
 bundle $E$ on $\Sigma$.
We then define by means of the Chern-Weil approach,
what should be the characteristic class of the would-be $d$-gerbe $\Det(q_*(E))$:
$$C_1(q_*(E)) \quad\in\quad H^{d+2}(B, \mathbb{C}).\leqno (0.4.1)$$
We denote it by $C_1$ since it is a kind of $d$-fold delooping of the usual first Chern
(determinantal) class. 
We then show the compatibility of this class with the gerbe approach whenever the latter
can be carried out rigorously. Our main result is the Real Riemann-Roch theorem (RRR):
$$C_1(q_*E) = \int_{\Sigma/B}\biggl[ ch(E)\cdot \Td({\mathcal T}_{\Sigma/B})\biggr]_{2d+2}
\quad\in\quad 
H^{d+2}(B, {\mathbb C}).
\leqno (0.4.2)$$
Here ${\mathcal T}_{\Sigma/B}$ is the complexified relative tangent bundle
and $\int_{\Sigma/B}$, the integration along the fibers of $q$, lowers the degree by $d$.

\vskip .2cm

 Note that the above theorem is a statement of purely real geometry and is
quite different from the ``Riemann-Roch theorem for differentiable manifolds'' proved by
Atiyah and Hirzebruch [AH]. The latter expresses properties of a Dirac operator on a real manifold
$\Sigma$, while our RRR deals with the $\overline{\partial}$-operator on a
complex envelope $X$ of $\Sigma$. The $d=1$ case above  can be deduced from  a result of Lott [Lo]
on ``higher" index forms for Dirac operators (because the polarization in the circle case
can be described in terms of the signs of eigenvalues of the Dirac operator). In general,
however, our results proceed in a different direction. 

\vskip .3cm

\noindent {\bf (0.5)} 
Our definition of $C_1(q_*E)$ uses the description of the cyclic homology of differential operators
[BG] [W] which provides a construction of a natural Lie algebra cohomology class
$\gamma$  of the Atiyah algebra,
i.e., of the Lie algebra of infinitesimal automorphisms of a pair
$(\Sigma, E)$ where $\Sigma$ is a compact oriented $d$-dimensional $C^\infty$-manifold
and $E$ is a vector bundle on $\Sigma$. The intuition with higher gerbes suggests that this class
comes in fact from a group cohomology class of the infinite-dimensional group of all the automorphisms
of $(\Sigma, E)$, see (3.7.7) and, moreover, that there
are similar classes coming from the higher Chern classes (3.7.8). 
 This provides a new point of view on the rather classical subject of
``cocycles on gauge groups and Lie algebras" i.e., 
on groups of diffeomorphisms of manifolds and
automorphisms of vector bundles as well as their Lie algebra analogs.

 There have been two spurs of interest in this subject.
The first one was the study of the cohomology of the Lie algebras of vector fields following
the work of Gelfand-Fuks, see [F] for a systematic account. In particular, Bott [Bo]
produced a series of cohomology classes of the Lie algebra of vector fields on a compact manifold
and integrated them to group cohomology classes of the group of diffeomorphisms. Later,
group cocycles have been studied with connections with various anomalies in physics,
see [RSF]. 
   
   From our point of view, the approach of
[RSF]  can be seen as producing  ``integrals of products of Chern classes"
in families over a base $B$, (cf. [D1] [E]),  in other words, as producing the
ingredients for the right hand side of a group-theoretical RRR.
This is the same approach that leads to the construction of the Morita-Miller
characteristic classes for surface fibrations [Mo]. 
The anomalies themselves, however,  should be seen as the classes whose existence is conjectured
in (3.7.7-8) and whose description through integrals of products of Chern classes constitutes
the  RRR. 

\vskip .3cm

\noindent {\bf (0.6)}  As far as the proof of the RRR goes, we use two types of techniques.
The first is that of differential graded Lie algebroids (which can be seen as infinitesimal
analogs of higher groupoids appearing in the heuristic discussion above). The second technique is
that of ``formal geometry" of Gelfand and Kazhdan, i.e., reduction of global problems in geometry
of manifolds  and vector bundles to problems related to cohomology of Lie algebras of  formal
vector fields and currents. The first work relating Riemann-Roch to Lie algebra cohomology was [FT]
and this approach was further developed in [BNT]. To prove the RRR we use results of 
[NT] and [BNT] on the Lie algebra
cohomology of formal Atiyah algebras.

\vskip .3cm

\noindent {\bf (0.7)} We are grateful to K. C. H. Mackenzie for 
pointing out several inaccuracies in an earlier version.
 The second author would like to acknowledge support from NSF,
Universit\'e Paris-7 and Max-Planck Institut fuer Mathematik.

\vfill\eject

\centerline {\bf 1. Background on Lie algebroids, groupoids and gerbes.}

\vskip 1cm

\noindent {\bf (1.1) Conventions.}
All manifolds will be understood to be $C^\infty$ unless otherwise specified. For a
manifold $\Sigma$ we denote by $C^\infty_\Sigma$ the sheaf of $\mathbb C$-valued
$C^\infty$-functions. By a vector bundle over $\Sigma$ we mean a locally trivial,
$C^\infty$ complex  vector bundle, possibly infinite-dimensional. For such a bundle
$E$ we denote by $C^\infty(E)=C^\infty_\Sigma(E)$ the sheaf of smooth sections,
which is a locally free sheaf of $C^\infty_\Sigma$-modules. By ${\mathcal T}_\Sigma$
we denote the {\em complexified} tangent bundle of $\Sigma$, so its sections are
derivations of $C^\infty_\Sigma$. We denote by $\cD_\Sigma$ the sheaf of differential
operators acting on $C^\infty_\Sigma$ and by $\cD_{\Sigma, E}$ the sheaf of differential
operators acting from sections of $E$ to sections of $E$. The notations $\cD(\Sigma)$ and
$\cD(\Sigma, E)$ will be used for the spaces of global sections of
$\cD_\Sigma$ and $\cD_{\Sigma, E}$. 

\vskip .3cm

\noindent {\bf (1.2) Lie algebroids.} Recall [Mac] that a Lie algebroid on $\Sigma$
consists of a vector bundle $\mathcal G$, a morphism of vector bundles $\alpha:{\mathcal G}
\to{\mathcal T}_\Sigma$ (the anchor map) and a Lie algebra structure in $C^\infty({\mathcal G})$
satisfying the properties: 

\vskip .2cm

\noindent (1.2.1) $\alpha$ takes the Lie bracket on sections of
$\mathcal G$ to the standard Lie bracket
on vector fields.

\vskip .1cm

\noindent (1.2.2) For any smooth function $f$ on $\Sigma$ and sections $x,y$ of $\mathcal G$
we have
$$[fx,y] - f\cdot [x,y] = \operatorname{Lie}_{\alpha(y)}(f)\cdot x.$$

\vskip .1cm

A Lie algebroid is called transitive, if $\alpha$ is surjective. 

\vskip .2cm

\noindent {\bf (1.2.3) Examples.} (a) When $\Sigma=pt$, a Lie algebroid is the same
as a Lie algebra. 

\vskip .1cm

(b) ${\mathcal T}_X$ with the standard Lie bracket and $\alpha=\id$  is a Lie
algebroid.

\vskip .1cm

(c) If $\alpha=0$, then the bracket in $\mathcal G$ is $C^\infty_\Sigma$-linear.
In this case we say that $\mathcal G$ is a bundle of Lie algebras: every fiber of
$\mathcal G$ is a Lie algebra. 

\vskip .2cm

For a fixed $\Sigma$ we will speak about morphims of Lie algebroids on $\Sigma$,
understanding morphisms of Lie algebroids in the sense of [Mac] which are identical
on $\Sigma$. Thus a morphism $\cG\to \cG'$ is a morphism of vector bundles
commuting with brackets and the anchor maps.

 Note that for any transitive Lie
algebroid $\mathcal G$ the kernel $\Ker(\alpha)\subset {\mathcal G}$ is a bundle
of Lie algebras, i.e.,  a Lie algebroid with trivial anchor map, and the maps in the short
exact sequence
$$0\to \Ker(\alpha)\to {\mathcal G}\buildrel \alpha\over\to {\mathcal T}_X \to 0
\leqno (1.2.4)$$
are morphisms of Lie algebroids.

\vskip .3cm

\noindent {\bf (1.3) The de Rham complex of a Lie algebroid.} Let
$\mathcal G$ be a Lie algebroid on $\Sigma$. We denote 
$$\DR^i({\mathcal G}) = \Hom(\Lambda^i{\mathcal G}, C^\infty_\Sigma). $$
The differential $d: \DR^i({\mathcal G}) \to \DR^{i+1}({\mathcal G})$
is defined by the standard formula of Cartan: for an antisymmetric $i$-linear
function $l: {\mathcal G}^i\to C^\infty_\Sigma $ we set
$$dl(x_1, ..., x_{i+1}) = \sum_{j=1}^{i+1} (-1)^j l(x_1, ...,\widehat{x_j}, ... , x_{i+1}) + 
\sum_{j<k} (-1)^{j+k} l([x_j, x_k], x_1, ..., \widehat{x_j}, ..., \widehat{x_k}, ..., 
x_{i+1}). \leqno (1.3.1).$$
We get a complex $\DR^\bullet({\mathcal G})$ called the de Rham complex of
$\mathcal G$.  A morphism of Lie algebroids $\phi: {\mathcal G}\to  
{\mathcal H}$ gives rise to the morphism of de Rham complexes
$\phi^*: \DR^\bullet (\cH)\to\DR^\bullet (\cG)$.

\vskip .3cm

\noindent {\bf (1.3.2) Examples.} (a) if $\Sigma=pt$, so $\cG$ is a Lie algebra,
then $\DR^\bullet(\cG) =  C^\bullet(\cG)$ is the cochain complex of
$\cG$ with trivial coefficients. 

\vskip .1cm

(b) If $\cG = \cT_\Sigma$, then $\DR^\bullet(\cG) = \Omega^\bullet_\Sigma$
is the $C^\infty$ de Rham complex of $\Sigma$. 

\vskip .3cm

\noindent {\bf (1.4) The enveloping algebra of a Lie algebroid.}
Let $\cG$ be a Lie algebroid on $\Sigma$, as before. The enveloping algebra
$U(\cG)$ is the sheaf of associative algebras on $\Sigma$ defined by generators
$x\in\cG$ (local sections) and $f\in C^\infty_\Sigma$ (local functions) subject to
the relations:
$$xy-yx = [x,y]. \leqno (1.4.1)$$
$$f\cdot x - x\cdot f = \Lie_{\alpha(x)}(f).\leqno (1.4.2)$$

\vskip .2cm

\noindent {\bf (1.4.3) Examples.} (a) If $\Sigma=pt$, so $\cG$
is a Lie algebra, then $U(\cG)$ is the usual enveloping algebra of $\cG$.

\vskip .1cm

(b) If $\cG = \cT_\Sigma$, then $U(\cG) = \cD_\Sigma$ is the sheaf of differential
operators $C^\infty_\Sigma\to C^\infty_\Sigma$.

\vskip .1cm

(c) If $\cG$ is any Lie algebroid, then the anchor map $\alpha$ induces a morphism
$$U(\alpha): U(\cG) \to U(\cT_\Sigma) = \cD_\Sigma$$
of sheaves of associative algebras. In particular, $C^\infty_\Sigma$ is a left
$U(\cG)$-module. 

\vskip .2cm

The sheaf $U(\cG)$ has an increasing ring filtration $\{U^m(\cG)\}$ with $U^m(\cG)$ generated
by products involving at most $m$ sections of $\cG$. The following is then standard.

\proclaim (1.4.4) Proposition. The associated graded sheaf of algebras $\operatorname{gr} \, U(\cG)$
is identified with the symmetric algebra $S^\bullet(\cG)$.

\vskip .3cm

\noindent {\bf (1.5) The Koszul resolution.} Let $\cG$ be a Lie algebroid on $\Sigma$.
We have then the complex
$$...\to U(\cG)\otimes \Lambda^2\cG \to U(\cG)\otimes\cG \to U(\cG) \to C^\infty_\Sigma\to 0.
\leqno (1.5.1)$$
with the differential defined by: 
$$d(u\otimes (\gamma_1\wedge ...\wedge \gamma_n)) = \sum_{j=1}^n (-1)^j (u\gamma_j) \otimes \bigl(\gamma_1\wedge ... \wedge
\widehat{\gamma_j}\wedge ... \wedge \gamma_n\bigr) + $$
$$+ \sum_{j<k} (-1)^{j+k} u\otimes \bigl([\gamma_i, \gamma_j] \wedge ... \wedge \widehat{\gamma_i} \wedge ...
\wedge \widehat{\gamma_j} \wedge ... \wedge \gamma_n\bigr).$$

\vskip 3cm

\proclaim (1.5.2) Proposition. The complex (1.5.1) is exact and thus provides a locally free
resolution of $C^\infty_\Sigma$ as a $U(\cG)$-module.

\proclaim (1.5.3) Corollary. We have
$$\DR^\bullet (\cG) \simeq \underline{R\Hom}_{U(\cG)} (C^\infty_\Sigma, 
C^\infty_\Sigma). $$

\vskip .3cm

\noindent {\bf (1.6) The Atiyah algebra.}
Let $G$ be a Lie group, $\Gg$ be its Lie algebra, and $\rho: P\to\Sigma$ a principal
$G$-bundle on $\Sigma$. The Atiyah algebra $\cA_P$ is the sheaf of Lie algebras on
$\Sigma$ whose sections are $G$-invariant vector fields on $P$:
$$\cA_P= (\rho_* \cT_P)^G. \leqno (1.6.1)$$
The map $\alpha=d\rho$ makes $\cA_P$ into a transitive Lie algebroid of the form
$$0\to\Ad(P) \to \cA_P\buildrel\alpha\over \longrightarrow \cT_\Sigma\to 0.\leqno (1.6.2)$$
Here $\Ad(P)$ is the bundle of Lie algebras on $\Sigma$ associated to $P$ via the adjoint
representation. 

If $\Sigma = \bigcup U_i$ is a covering in which $P$ is trivialized:
$P|_{U_i} = U_i\times G$, and $g_{ij}: U_i\cap U_j\to G$
are the transition functions, then $\cA_P$ is glued out of
$\cA_P|_{U_i} = \cT_{U_i}\times \Gg$ via the transition functions 
$$(v,x) \mapsto (v,  i_v(dg_{ij}\cdot g_{ij}^{-1}) + \Ad_{g_{ij}}(x)).\leqno (1.6.3)$$

\vskip .1cm

\noindent {\bf (1.6.4) Example.} Let $G= GL_r(\bbC)$, so $\Gg = {\frak {gl}}_r(\bbC)$. 
A principal $G$-bundle $P$ corresponds then to a rank $r$ vector bundle $E$ on $\Sigma$. In this
case $\cA_P$ will also be denoted $\cA_E$ and has a well known alternative description. It consists of 
differential operators $L: E\to E$ such that:

\vskip .2cm

\noindent (a) $L$ has order $\leq 1$.

\vskip .1cm

\noindent (b) The first order symbol of $L$ (which is a priori a section of
$\cT_\Sigma\otimes \End(E)$) lies in the subsheaf $\cT_\Sigma = \cT_\Sigma\otimes 1$.

\vskip .3cm

\noindent {\bf (1.7) Modules over  Lie algebroids.} 
We follow [Mac], see also [Kal] \S 3 for  a more algebraic language. 
Let $\cG$ be a Lie algebroid on 
$\Sigma$. A $\cG$-module is a vector bundle $\cM$ on $\Sigma$ equipped with a Lie algebra action
$(x, m)\mapsto xm$ of $\cG$ in the sections which satisfies the two twisted linearity properties:
$$ x(f\cdot  m) - f\cdot (xm) = (\Lie_{\alpha(x)} f)\cdot m, \quad
f\in C^\infty_\Sigma, x\in\cG, m\in\cM.\leqno (1.7.1)$$

\noindent (1.7.2) For any local section $x$ of $\cG$ the operator $m\mapsto xm$ on sections of
$\cM$ belongs to the Atiyah algebra $\cA_\cM$, and for different $m$ these operators define a
morphism of Lie algebroids $\cG\to\cA_\cM$.

\vskip .2cm

\noindent {\bf (1.7.3) Examples.} (a) For any $\cG$ the trivial bundle (whose sheaf of sections is)
$C^\infty_\Sigma$ is a $\cG$-module with the $\cG$ action given via the anchor map and the Lie derivations
of functions.

\vskip .1cm

(b) Unless the anchor of $\cG$ is trivial, the bracket in $\cG$ does not make $\cG$ into a
$\cG$-module, as the map from $\cG$ to the  Atiyah algebroid of $\cG$ as a vector bundle
is not linear over functions.

\vskip .1cm

(c) An ideal in $\cG$  is a sub-Lie algebroid $\cG'$ such that $[\cG, \cG']\subset \cG'$.
In this case $\cG'$ is a $\cG$-module provided its anchor is trivial. 

\vskip .2cm

Any $\cG$-module has a structure of a sheaf of modules over the sheaf of rings $U(\cG)$.

\vskip .3cm

\noindent {\bf (1.8) Cohomology of Lie algebroids.}  
Let $\cM$ be a $\cG$-module. The
  de Rham complex $\DR^\bullet(\cG, \cM)$ with coefficients in $\cM$ is defined by
$$\DR^i(\cG, \cM) = \underline {\Hom}(\Lambda^i{\mathcal G}, \cM). \leqno (1.8.1)$$
with the differential of $l: \cG^i\to \cM$ defined by the modification of (1.3.1): 
$$dl(x_1, ..., x_{i+1}) = \sum_{j=1}^{i+1} (-1)^j x_j (l(x_1, ...,\widehat{x_j}, ... , x_{i+1})) + 
\sum_{j<k} (-1)^{j+k} l([x_j, x_k], x_1, ..., \widehat{x_j}, ..., \widehat{x_k}, ..., 
x_{i+1}). \leqno (1.8.2).$$
Its cohomology sheaves will be denoted $\underline{H}_{\Lie}^i(\cG, \cM)$
and the corresponding cohomology groups of the complex of global smooth sections of $\DR^\bullet(\cG, \cM)$ by simply 
$H^i_{\Lie}(\cG, \cM)$. See [Mac], \S 7.1. 
As before, it is easy to see that
$$\DR^\bullet(\cG, \cM) \simeq \underline{R\Hom}_{U(\cG)} (C^\infty_\Sigma, \cM).\leqno (1.8.3)$$
Therefore 
$$\underline{H^i}_{\Lie}(\cG, \cM) = \underline{\Ext}^i_{U(\cG)}(C^\infty_\Sigma, \cM),\quad
H^i_{\Lie}(\cG,\cM) = \Ext^i_{U(\cG)}(C^\infty_\Sigma, \cM).\leqno (1.8.4)$$

\vskip .1cm

\noindent {\bf (1.8.5) Example.} 
The trivial bundle  $C^\infty_\Sigma$ is always a $\cG$-module and for $\cG = \cT_\Sigma$
we have $H^i_{\Lie}(\cT_\Sigma, C^\infty_\Sigma) = H^i(\Sigma, \bbC)$
(topological cohomology). 

\vskip .3cm

\noindent {\bf (1.9) The Hochschild-Serre spectral sequence and the transgression.} 
Let
$$0\to\cG'\to\cG\to\cG''\to 0\leqno (1.9.1)$$
be an extension of Lie algebroids on $\Sigma$, so $\cG'$ is an ideal in $\cG$.
Note that $\cG'$ is then a bundle of Lie algebras. Let $\cM$ be a $\cG$-module.
Then for every point $x\in\Sigma$ the fiber $\cM_x$ is a module over the
Lie algebra $\cG'_x$. Assume that for any $i\geq 0$ the Lie algebra
cohomology spaces
$H^i_{\Lie}(\cG'_x, \cM_x)$ have dimension independent on $x\in\Sigma$.
Then the sheaves $\underline{H}^i_{\Lie} (\cG', \cM)$ are vector bundles on
$\Sigma$ and these vector bundles have natural structures of
$\cG''$-modules. In this case
we have (a Lie algebroid generalization of)
 the Hochshild-Serre spectral sequence with
$$E_2^{pq} = H^p_{\Lie}(\cG'', \underline{H}^q_{\Lie}(\cG', \cM))\Rightarrow
H^{p+q}_{\Lie}(\cG, \cM).\leqno (1.9.2)$$
The construction is parallel to the classical (Lie algebra) case as in [F]. One uses
the short exact sequence (1.9.1) to produce, in a standard way,  a filtration on
$\DR^\bullet(\cG, \cM)$. See [Mac], \S 7.4 for the treatment of the case $\cG'' = \cT_\Sigma$
which is the only case we will use in this paper. 

\vskip .1cm

\noindent {\bf (1.9.3) Example.}
Similarly to the classical case, one can  use (1.9.2) (or elementary considerations)
to identify $H^2_{\Lie}(\cG, \cM)$ with the
set of isomorphism classes of central extensions of Lie algebroids
$$0\to \cM\to\widetilde{\cG}\to\cG\to 0.$$
Central extensions of this type with  $\cG = \cT_\Sigma$,  $\cM= C^\infty_\Sigma$,
and the $\cG$-action on $\cM$ being the standard one (by Lie derivations),
were called in [Kal] Picard Lie algebroids. The set of their isomorphism classes
is thus identified with  $H^2_{\Lie}(\cT_\Sigma, C^\infty_\Sigma)$ which is the same as
 the topological (de Rham) 
cohomology $H^2(\Sigma, \bbC)$. 

\vskip .1cm
Fix $n>0$ and assume that 
$$H^j(\cG', \cM) = 0, \quad 0<j<n, \leqno (1.9.4)$$
In this case $E_2^{0,n}=E_{n+1}^{0,n}$ as well as
$E_2^{0, n+1} = E_{n+1}^{0, n+1}$. We obtain therefore the
{\em transgression map}
$$d_{n+1}: E_{n+1}^{0,n} = E_2^{0,n} = 
H^n_{\Lie}(\cG', \cM)^{\cG''} \to H^{n+1}_{\Lie}(\cG'', \cM^{\cG'})= E_2^{n+1, 0} = E_{n+1}^{n+1, 0}.
\leqno (1.9.5)$$
We will use this map later in the paper. Without the assumption (1.9.4) we have 
 that $E_{n+1}^{0,n}$ is a subspace of 
$E_2^{0,n} = H^n_{\Lie}(\cG', \cM)^{\cG''}$ 
namely the intersection of the kernels of $d_2, ..., d_n$. For convenience we will
call elements of this space {\em transgressive}  elements of $E_2^{0,n}$.  
Similarly, $E_{n+1} ^{n+1, 0}$ is a quotient space of
$E_2^{n+1, 0} = H^{n+1}_{\Lie}(\cG'', \cM^{\cG'})$  by the
union of images of $d_2, ..., d_n$.  

\vskip .2cm

\noindent {\bf (1.9.6) Example.} Suppose that $n=2$ and $\Sigma = pt$, so (1.9.1)
is a central extension of Lie algebras and $\cM$ is a $\cG$-module in the
usual sense. Let $\gamma\in E_2^{0,2} = H^2_{\Lie}(\cG', \cM)^{\cG''}$ be
a $\cG''$-invariant class in $H^2$ and 
$$0\to\cM\to \widetilde{\cG}' \to \cG' \to 0$$
be a central extension representing $\gamma$. The class $\gamma$ is transgressive,
(i.e., annihilated by $d_2$) if and only if $\widetilde{\cG'}$ can be made into
a $\cG$-equivariant central extension
 (as opposed to the fact that the class of the extension remains unchanged under the
$\cG$-action or, what is the same, under $\cG''$-action). 
 Given such an equivariant extension, one obtains a crossed module of Lie algebras
 (i.e., a dg-Lie algebra situated in degrees (-1) and 0)
$$ \widetilde{\cG}''\buildrel\partial\over\longrightarrow \cG,$$
with $\Ker(\partial) = \cM$ and $\Coker(\partial) = \cG''$. As well known
(see, e.g., [L], Example E.10.3), 
such a crossed module represents an element in $H^3(\cG'', \cM)$, and this element
is the lifting of  $d_3(\beta)$.
Different choices of equivariant structure on $\widetilde{\cG'}$ 
correspond to the ambiguity of the values of $d_3$ modulo the image of $d_2$. 
One can generalize this picture easily to the case of an arbitrary $\Sigma$. 

\vskip .3cm

\noindent {\bf (1.10) Reminder on gerbes.} We follow the same conventions as in [KV2] and
use [Bre] as the background reference.

If $B$ is a topological space and  
and $\cF$ is a sheaf of abelian groups on $B$, then we can speak of 
$\cF$-gerbes (= gerbes with lien $\cF$). Recall that such a 
gerbe ${\frak{ G}}$ 
consists of the following data: 

\vskip .2cm

(1) A category ${\frak{ G}}(U)$ given for all open $U\subset B$, 
the restriction functors $r_{UV}: {\frak{G}}(U)\to {\frak {G}}(V)$ given for any 
morphism 
$V\subset U$ 
and natural isomorphisms of functors $s_{UVW}: r_{VW}\circ r_{UV}\Rightarrow 
r_{UW}$ 
given for each $W\subset V\subset  U$ and satisfying the transitivity conditions. 

\vskip .1cm

(2) The structure of ${\cF}|_U$-torsor (possibly empty) on each sheaf
$\shHom_{{\frak {G}}(U)} 
(x,y)$ compatible with the $r_{UV}$ and such that the composition of morphisms 
is bi-additive. 

\vskip .2cm

These data have to satisfy the local uniqueness and gluing properties for which 
we refer to [Bre]. 

By a sheaf of $\cF$-groupoids we will mean a sheaf of categories ${\frak{ C}}$ 
on $B$ 
(so both ${ \Ob}\, {\frak {C}}$ and 
$\Mor \, {\frak {C}}$ are sheaves of sets) 
in which each sheaf $\shHom_{{\frak {C}}(U)}(x,y)$ is either empty or is made 
into a sheaf of  ${\cF}|_U$-torsors 
so that the composition is biadditive. 
A sheaf ${\frak {C}}$  of $\cF$-groupoids is called 
locally connected if locally on $B$ all the ${\Ob} \, {\frak {C}}(U)$ 
and  $\Hom_{{\frak {C}}(U)}(x,y)$ 
are nonempty.

Each sheaf of $\cF$-groupoids can be seen as a fibered category 
over $B$, in fact it is a pre-stack, see, e.g., [LM]. 
Recall (see, e.g., {\it loc. cit.} Lemma 2.2) that for any pre-stack 
${\frak {C}}$ there is an associated stack ${\frak {C}}^{^\sim}$. 
If ${\frak {C}}$ is a locally connected sheaf of $\cF$-groupoids, 
then ${\frak {C}}^{^\sim}$ is an $\cF$-gerbe.   

As well known (see, e.g., [Bre]), the set formed by 
${\cF}$-gerbes up to equivalence is identified with $H^2(B, {\cF})$. 
The identification of the set of isomorphism classes of Picard Lie
algebroids in Example 1.9.3 can be seen as an infinitesimal analog of this fact. 
Given an $\cF$-gerbe ${\frak{ G}}$, we denote by $[{\frak {G}}]\in H^2(B, {\cF})$ 
its class. Given a sheaf ${\frak{ C}}$ of $\cF$-groupoids, we denote 
by $[{\frak {C}}]$ the class of the corresponding gerbe. 

\vskip .2cm

Let $B$ be a $C^\infty$-manifold. We will be particularly interested in $C^{\infty *}_B$-gerbes
on $B$. 
Recall that we have the exponential sequence of sheaves on $B$:
$$0\to \underline {\bbZ}_B \to  C^\infty_B\buildrel e^{2\pi i x} \over\longrightarrow
C^{\infty *}_B\to 0.\leqno (1.10.1)$$
The corresponding coboundary map
$$\delta_n: H^n(B, C^{\infty *}_B\to H^{n+1}(B, \bbZ)\leqno (1.10.2)$$
is an isomoprhism for $n\geq 1$ since $C^\infty_B$ is a soft sheaf. 
Thus $[{\frak {G}}]$ give rise to a class in $H^3(B, \bbZ)$. 

\vskip .2cm

Let ${\frak {G}}$ be a $C^{\infty *}_B$-gerbe.
Recall [Bry], that a {\em connective structure}  $\Delta$ on ${\frak {G}}$
is a set of data that associates to each open $U\subset B$ and each object
$x\in \Ob \, {\frak{G}}(U)$ a sheaf $\Delta(x)$ of  $\Omega^1_U$-torsors
(whose sections can be thought of as ``formal connections" in $x$) and for any
local (iso)morphism $g: x\to y$ over $U$ an identification of torsors
$g_*: \Delta(x)\to\Delta(y)$, satisfying the compatibility property plus the 
following gauge condition: if $x=y$ so $g\in C^{\infty *}(U)$ is an invertible
function, then $g_*(\nabla) = \nabla - g^{-1} d(g)$.

A {\em curving} of a connective structure $\Delta$ is a rule $K$ associating
to any $x$ as above and any global object $\nabla\in\Delta(x)$ a 2-form
$K(\nabla)\in\Omega^2(U)$ satisfying the compatibility with pullbacks,
invariance under isomorphisms as well as the gauge condition:
$K(\nabla + \alpha) = K(\nabla) + d\alpha$, $\alpha\in\Omega^1(U)$.
 In this situation Brylinski defined the 3-curvature of the
connective structure and curving, which 
is a closed 3-form $S = S_{\Delta, K}\in\Omega^3(B)$. 

\vskip .2cm

\noindent {\bf (1.10.3) Example.} let $G$ be a Lie group and
$$1\to \bbC^*\to \widetilde{G}\to G\to 1$$
be a central extension of Lie groups. Let $\rho: P\to B$ be a principal
$G$-bundle. We have then the $C^{\infty *}_B$-gerbe $\Lift_G^{\widetilde{G}}(P)$
whose  objects  over $U\subset B$  are liftings of $P|_U$ to a principal 
$\widetilde{G}$-bundle over $U$, compare [Bl]. Let $\nabla_P$ be a connection on
$P$. Then $\Lift_G^{\widetilde G}(P)$ has a connective structure $\Delta$
which to every lifting $\widetilde {P}$ of $P$ to a $\widetilde{G}$-bundle
associates the space of all connections on $\widetilde{P}$ extending $\nabla_P$.
Further, let $R_\nabla\in\Omega^2(B)\otimes \Ad(P)$ be the curvature of $\nabla$.
A choice of a lifting of $R_\nabla$ to a form $\widetilde{R}_\nabla\in\Omega^2(B)\otimes
\Ad(\widetilde{P})$ gives a curving $K$ on $\Delta$. This curving associates
to any section $\widetilde{\nabla}$ of $\Delta(\widetilde{P})$, i.e., to
a connection on $\widetilde{P}$ extending $\nabla$, the 2-form $R_{\widetilde {\nabla}} - 
\widetilde{R}_\nabla$, where $R_{\widetilde{\nabla}}$ is the curvature of 
$\widetilde{\nabla}$.

\vskip .2cm

We will need the following result ([Bry], Thm. 5.3.12).

\proclaim (1.10.4) Theorem. If ${\frak{G}}$ is a $C^{\infty *}_B$-gerbe
with a connective structure $\Delta$ and a curving $K$, then the class of 
 $S_{\Delta, K}$ in $H^3(B, \bbC)$ is integral and is equal to the image
 of $\delta_2 [{\frak {G}}]$ under the natural map from $H^3(B, \bbZ)$
 to $H^3(B, \bbC)$.

\vfill\eject

\centerline {\bf 2. Background on homology of differential operators. }

\vskip 1cm

\noindent{\bf (2.1) Conventions.} Let $A$ be an associative algebra over $\CC$. We denote by
$\Hoch_\bullet(A)$ the Hochschild complex of $A$ with coefficients in $A$:
$$...\rightarrow A\otimes A\otimes A \to A\otimes A \to A,\leqno (2.1.1)$$
$$d(a_0\otimes ...\otimes a_p) = \sum_{i=0}^{p-1} (-1)^i a_0\otimes ... \otimes a_i a_{i+1}\otimes ... \otimes a_p
+
(-1)^p a_p a_0 \otimes a_1\otimes ... \otimes a_{p-1}.$$
By $HH_\bullet (A)$ wew denote the homology of $\Hoch_\bullet(A)$. As well known,
$$HH_\bullet(A) = \operatorname{Tor}^{A\otimes A^{op}}_\bullet(A, A).\leqno (2.1.2)$$
Put
$$\tau(a_0\otimes ... \otimes a_p) = (-1)^p a_1\otimes ... \otimes a_p\otimes a_0.\leqno (2.1.3)$$
The cyclic complex of $A$ is defined as the total compplex
$$CC_\bullet(A) = \operatorname{Tot}_\bullet \biggl\{ \cdots\Hoch_\bullet(A)\buildrel 1-\tau\over\rightarrow
\Hoch_\bullet(A)\buildrel N\over\rightarrow \Hoch_\bullet(A)\buildrel 1-\tau\over\rightarrow
\Hoch_\bullet(A)\biggr\}.
 \leqno (2.1.4)$$
Here   $N= 1+\tau+\tau^2+...+\tau^n$ on $\Hoch_n(A)$. 
 
The cyclic homology  $HC_\bullet(A)$ is the homology of the complex $CC_\bullet(A)$.
We recall the fundamental result relating the cyclic homology with the Lie
algebra homology of the algebra of matrices, see [L].

\proclaim (2.1.5) Theorem.
$$H_\bullet^{Lie} ({\frak gl}(A)) = S^\bullet(HC_{\bullet-1}(A)).$$

\proclaim (2.1.6) Corollary.  If $HC_j(A)=0$ for $j=0, ..., p-1$, then $H_j^{Lie}({\frak gl}(A))=0$
for $j=0, ..., p$, and $H_{p+1}^{Lie}({\frak gl}(A)) = HC_p(A)$.

\vskip .3cm

\noindent {\bf (2.2) Homology of differential operators: algebro-geometric version.}
Let $X$ be a smooth affine algebraic variety over $\CC$ of dimension $d$ and $\cE$
be an algebraic vector bundle on $X$. Then the Hochschild-Kostant-Rosenberg theorem
(together with Morita invariance of $HH_\bullet$) gives an identification:
$$HH_m(\End(\cE)) = \Omega^m(X), \leqno (2.2.1)$$
where on the right we have the space of global regular $m$-forms on $X$.
 Further,
$$HC_m(\End(\cE)) = \Omega^m(X) d\Omega^{m-1}(X)
 \oplus H^{m-2}(X, \bbC ) \oplus H^{m-4}(X, \bbC) \oplus ..
\leqno (2.2.2)$$
where on the right we have the usual topological (de Rham) cohomology,
see [L] Th. 3.4.12. 
Let $\cD(\cE)$ be the ring of global differential operators from $\cE$ to $\cE$. Then the results
of [BG] [W] imply:
$$HH_m(\cD(\cE))= H^{2d-m}(X, \CC).\leqno (2.2.3)$$
 Further,
$$HC_m(\cD(\cE)) = \bigoplus_{i\geq 0} H^{2d-m+2i}(X, \CC). \leqno (2.2.4)$$
We recall that the approach of {\em loc. cit} is to use the filtration by the degree of differential
operators and realize the $E_1$-term of the corresponding spectral sequence for $HH$
as the complex of form on the cotangent bundle with the differential adjoint to the
de Rham differential by means of the symplectic form. The spectral sequence is then
seen to degenerate at $E_2$.

Let us note the particular case when $X=\bbA^d$ and $E= \cO_{\bbA^d}$
is the trivial bundle of rank 1. Then $\cD(\cE) = W_d$ is the Weyl algebra with generators
$x_i, \partial_i$, $i=1, ..., d$, and relations
$$[x_i, x_j] = [\partial_i, \partial_j] = 0, \quad [\partial_i, x_j] = \delta_{ij}\cdot 1.$$
The above results imply that

$$HH_i(W_d) =  0 \quad \text {if} \quad i\neq 2d, \quad \quad HH_{2d}(W_d) = \CC. \leqno (2.2.5)$$
and
$$HC_{i}(W_d) = \CC, \, i-2d\in 2\bbZ_+, \quad  HC_i(W_d)=0, \, i-2d\notin 2\bbZ_+ .\leqno (2.2.6)$$ 

\vskip .2cm

\noindent {\bf (2.3) The $C^\infty$ version.} Let $\Sigma$ be an oriented $C^\infty$-manifold of
dimension $d$ and $E$ be a smooth complex vector bundle on $\Sigma$. We have then the algebras
$\End(E)$, $\cD(E)$ of smooth endomorphisms and differential operators on $E$. Following
[W] we present the analogs of the results cited in (2.2) for these algebras.
 These rings have natural Fr\'ech\'et topologies. 
 As pointed out in loc. cit., to get reasonable results, all tensor
products occuring in the Hochschild and cyclic complexes of the above algebras
 should be taken taken in the category
of topological vector spaces, i.e., be completed.
In plain terms, this means that the $\End(E)^{\otimes p}$ should be understood as the ring of
endomorphisms of the vector bundle $E^{\boxtimes p}$ on the $p$-fold Cartesian product $\Sigma^p$
and similarly for differential operators. Under these conventions, we have:
$$HH_m(\cD(E))= H^{2d-m}(\Sigma, \CC), \leqno (2.3.1)$$
$$HC_m(\cD(E)) = \bigoplus_{i\geq 0} H^{2d-m+2i}(\Sigma, \CC), \leqno (2.3.2)$$
where on the right we have the topological cohomology. 

\vskip .2cm

\noindent {\bf (2.3.3) Remark.} The Lie algebra cochain complexes of $\cD(E)$ and of
${\frak {gl}}_N \cD(E) = \cD(E\otimes \bbC^r)$ involve exterior products of these
algebras over $\bbC$. If we understand these products in the completed sense as above
(compare also with Fuks [F]), then the analog of (2.1.5) holds, and we have the
following.

\proclaim (2.3.4) Corollary.  Let $\Sigma$ be a compact, oriented $C^\infty$ manifold
of dimension $d$. Then, for $N\gg 0$ we have:
$$H_i^{\Lie} {\frak{gl}}_N \cD(E)=0, \quad 0<i<d+1;$$
$$H_{d+1}^{\Lie} {\frak{gl}}_N\cD(E) = \bbC.$$

\vskip .2cm

\noindent {\bf (2.4) The formal series version.} Let
$$\widehat{W}_d =  W_d \otimes_{\CC[x_1, ..., x_d]} \CC[[x_1. ..., x_d]]$$
be the algebra of differential operators whose coefficients are formal power series.
Similarly to the 
 above, we consider the Hochschild and cyclic complexes of $\widehat{W}_d$
using the adic topology on $\CC[[x_1, ..., x_d]]$ and taking completions.
Thus $\widehat{W}_d^{\otimes p}$ is understood as the ring of differential operators
whose coefficients are power series in $p$ groups of $d$ variables. With this understanding,
we have the analog of (2.2.5):
$$HH_{2d}(\widehat{W}_d) = \CC, \quad HH_i(\widehat{W}_d)=0, \, i\neq 2d.\leqno (2.4.1)$$
For the proof, see  [FT]. One can also apply the spectral sequence argument of [BG] and [W]
and then use the Poincare lemma on the (contangent bundle to the) formal disk. 

\vskip .2cm

Our next step is to consider such formal completions simultaneously at all points of a
given $C^\infty$-manifold $\Sigma$. So let $\Sigma, E$ be as above. Let
$\widehat{\Hoch}_p (\cD(E))$ be the completion of $\cD(E^{\boxtimes(p+1)})$
(differential operators in the bundle $E^{\boxtimes (p+1)}$ on $\Sigma^{p+1}$)
along the diagonal $\Sigma\subset\Sigma^{p+1}$. This is a sheaf of $\Sigma$.

Then the Hochschild differential
extends to $\widehat{\Hoch}_\bullet(\cD(E))$, making it into 
a complex ,and we denote by 
its homology. Similarly, we define the completed cyclic complex $\widehat{CC}_\bullet(\cD(E))$
by the procedure identical to (2.1.4) and denote its homology by
$\widehat{HC}_\bullet(\cD(E))$. Thus, $\widehat{HH}_\bullet (\cD(E))$
and $\widehat{CC}_\bullet(\cD(E))$ are sheaves on $\Sigma$. 

\proclaim (2.4.2) Proposition. We have $\widehat{HH}_p(\cD(E)) =\underline{\CC}_\Sigma$
(constant sheaf) for $p=2d$
and it is equal to 0 for $p\neq 2d$. 

\noindent {\sl Proof:} It is clearly enough to consider the case when 
  $\Sigma$ is an open ball in $\bbR^d$
and $E$ is trivial and to prove that in this case the complex of global sections of
$\widehat{HH}_\bullet(\cD(E))$ is exact everywhere except degree $2d$ where the
cohomology is isomoprhic to $\bbC$. This will imply that the only sheaf of
cohomology of $\widehat{HH}_\bullet(\cD(E))$ in this case (and thus in the general case)
is $\underline{\bbC}_\Sigma$. So we make this assumption in the rest of the proof. 

We start with the case of $\widehat{Hoch}_\bullet(C^\infty_\Sigma)$ defined, as before,
using functions on the completion of $\Sigma^{\bullet +1}$ along the diagonal. 
Such functions form a flat module over $C^\infty(\Sigma\times\Sigma)$
so by the interperation of $HH$ as $\operatorname{Tor}$, see (2.1.2) we see that 
$$\widehat{HH}_\bullet(C^\infty_\Sigma) = \Omega^i_\Sigma,$$
and the same will hold if we replace $C^\infty_\Sigma$ by a matrix algebra
(i.e., take $E$ of higher rank).

Next, we replace $C^\infty_\Sigma$ by the sheaf of commutative algebras
$$\cA = S^\bullet(\cT_\Sigma)$$
(polynomial functions on the contangent bundle) and define $\widehat{Hoch}_\bullet(\cA)$
using the completions of sheaves of sections of $\cA^{\boxtimes (p+1)}$ on $\Sigma^{p+1}$
along the diagonals. The same argument with flatness  will apply, so we conclude that 
$$\widehat{HH}_\bullet(\cA) = p_* (\Omega^\bullet_{T^*\Sigma}),\leqno
(2.4.2.1) $$
where $p: T^*\Sigma\to \Sigma$ is the projection. Again, a similar statement will hold for
matrices. 

Finally, we use the approach of [BG] [W] and consider the spectral sequence for $\widehat{HH}_\bullet(\cD(E))$
associated to the filtration by degree of operators. We get the $E_1$-term to be (2.4.2.1)
with the differential being the adjoint of the de Rham differential on $T^*\Sigma$.
Since we assumed $\Sigma$ to eb a ball,
we conclude that the $E_2$-term reduces to one space $\bbC$. \qed

\vskip .2cm

Further, we need a relative version of the above statements. Let
$$q: \Sigma\to B$$
be a submersion (smooth fibration) of $C^\infty$-manifolds, whose fibers are of
dimension $d$ and are oriented. Let $E$ be $C^\infty$-bundle on $\Sigma$, as above.
We have then the subring
$$\cD_{\Sigma/B}(E) \subset \cD(E),\leqno (2.4.3)$$
consistsing of differential operators that are $q^{-1}C^\infty_B$-linear, i.e.,
act along the fibers only. 

Let $\Sigma^{p+1}_B\subset \Sigma^{p+1}$ be the $(p+1)$-fold fiber product of
$\Sigma$ over $B$. We denote by $E^{\boxtimes (p+1)}_B$ the restriction of
$E^{\boxtimes (p+1)}$ to $\Sigma^{p+1}_B$. 

Let $\widehat{\Hoch}_p(\cD_{\Sigma/B}(E))$ denote the completion of
$\cD_{\Sigma^{p+1}_B/B}(E^{\boxtimes (p+1)}_B)$ along the diagonal. Then the Hochschild
differential extends to $\widehat{\Hoch}_p(\cD_{\Sigma/B}(E))$. We also define the completed
cyclic complex $\widehat{CC}_\bullet(\cD_{\Sigma/B}(E))$ by implementing (2.1.4). 

\proclaim (2.4.4) Theorem. (a) The complex $\widehat{\Hoch}_p(\cD_{\Sigma/B}(E))$ is acyclic in
degrees other than $2d$, and its $2d$th cohomology sheaf is
isomorphic to $q^{-1} C^\infty_B$. In other words, we have a quasiisomorphism in the derived
category of sheaves of $q^{-1}C^\infty_B$-modules on $\Sigma$:
$$\mu_\cD: \widehat{\Hoch}_p(\cD_{\Sigma/B}(E)) \to q^{-1} C^\infty_B [2d].$$
(b) We have $H^i(\widehat{CC}_\bullet(\cD_{\Sigma/B}(E))=0$ unless $i=-2d+k$, $k\in \bbZ_+$, and
$$H^{-2d+k}(\widehat{CC}_\bullet (\cD_{\Sigma/B}(E)) = q^{-1} C^\infty_B.$$

\noindent {\sl Proof:} Similar to (2.4.2). 

\vskip .2cm

As a corollary of Theorem 2.4.4 we have a morphism (no longer an isomorphism) in the derived category
$$\nu_\cD: \widehat{CC}_\bullet(\cD_{\Sigma/B}(E))\to q^{-1} C^\infty_B[2d]. \leqno (2.4.5) $$

\vfill\eject

\centerline {\bf 3. Characteristic classes from Lie algebra cohomology.}

\vskip 1cm

\noindent {\bf (3.1) The finite-dimensional case.} Let $G$ be a 
Lie group with Lie algebra $\frak g$. We denote by $C^\bullet({\frak g})$
the cochain complex of $\frak g$ with trivial coefficients $\bbC$
and by $H^n({\frak g})$ its $n$th cohomology space. 

Let $\gamma\in H^n({\frak g})$ be a cohomology class. We want to associate (under certain
conditions) to $\gamma$ a characteristic class of principal $G$-bundles. In other words,
we want to produce,  for each $C^\infty$-manifold $B$ and each smooth principal $G$-bundle
$P$ on $B$, a topological (de Rham) cohomology class
$$c_\gamma(P)\in  H^{n+1}(B) = H^{n+1}(B, \bbC)  \leqno (3.1.1)$$
(note the shift of degree by $1$). 

\vskip .2cm

Indeed, let a principal $G$-bundle $\rho:  P\to B$ be given and $\cA_P$ be its Atiyah algebra.
We have then the extension of Lie algebroids (1.6.2) on $B$ and the corresponding Hochschild-Serre
spectral sequence (1.9.2) which in our case has the form:
$$E_2^{pq} = H^p_{\Lie}(\cT_B, \underline{H}^q_{\Lie}(\Ad(P), C^\infty_B))\Rightarrow H^{p+q}_{\Lie}(\cA_P, C^\infty_B).
\leqno (3.1.2)$$
This sequence was considered in [McK], Thm. 7.4.19. 
Note that $\underline{H}^q_{\Lie}(\Ad(P), C^\infty_B)$ is the
cohomology of the cochain complex of $\Ad(P)$ as a Lie algebra over $C^\infty_B$, 
 i.e., of
the complex of bundles formed by the duals of the fiberwise exterior products of fibers of $\Ad(P)$.
We will also use the notation $C^\bullet(\Ad(P)_{/B})$ for this complex.

\proclaim (3.1.3) Lemma. For any $q\geq 0$ the bundle $H^q_{\Lie}(\Ad(P), C^\infty_B) =
H^q(\Ad(P)_{/B})$ on $B$ formed by the 
Lie algebra cohomology spaces of the fibers of $\Ad(P)$ is canonically identified with the
trivial bundle with fiber $H^q(\Gg)$. 

\noindent{\sl Proof:} This follows from the fact the the adjoint action of $G$ on $\Gg$ induces
the trivial action on $H^q(\Gg)$.

\vskip .1cm

Therefore 
$$E_2^{pq} = H^p(B)\otimes H^{q}(\Gg). \leqno (3.1.4)$$
In particular, our class $\gamma\in H^n(\Gg)$ gives an element
$1\otimes\gamma\in E_2^{0n}$.



Assume now that we have $n>0$ such that
the Lie algebra $\Gg$ satisfies the sphericity condition:
$$H^i(\Gg) = 0, \quad 0<i<n. \leqno (3.1.5)$$
Then we are in the situation of (1.9.4), so we have the transgression map (1.9.5)
which in our case has the form
$$d_{n+1}: H^n(\Gg)\to H^{n+1}(B), \leqno (3.1.6)$$
and we define
$$c_\gamma(P) = d_{n+1}(1\otimes\gamma). \leqno (3.1.7)$$
Without the assumption (3.1.5) we have that $c_\gamma(P)$ is defined
only if $1\otimes\gamma$ is transgressive (i.e., annihilated by $d_2, ..., d_n$
and takes value not in $H^{n+1}(B)$ but in the quotient of $H^{n+1}(B)$
by the images of $d_2, ..., d_n$.

\vskip .2cm

\noindent {\bf (3.1.8) Examples.}  

 (a) Let $n=1$. Then the condition (3.1.5) is
trivially satisfied. A class $\gamma$ is just a trace functional $\gamma: \Gg\to\bbC$.
The class $c_\gamma(P)\in H^2(B)$ can be obtained by choosing a connection $\nabla$
in $P$ with curvature $R\in \Omega^2_B\otimes\Gg$ and taking the class of the closed
2-form $\gamma(R)\in\Omega^2_B$. Alternatively, one can use $\gamma$ to produce a trace functional
$\gamma_P: \Ad(P)\to C^\infty_B$ and then use $\gamma_P$ to push forward the extension
(1.6.2) to a central extension of Lie algebroids
$$0\to C^\infty_B\to \cG\to \cT_B\to 0.$$
As well known (1.7) the set of isomorphism classes of such central extensions is identified with
$H^2_{\Lie}(\cT_B, C^\infty_B) = H^2(B, \bbC)$. 

\vskip .1cm

(b) Let $n=2$, so $\gamma$ is represented by a central extension
$$0\to \bbC\to \widetilde{\Gg} \to \Gg \to 0.\leqno (3.1.9)$$
A sufficient condition for $\gamma$ to be basic for any $P$ is that $\widetilde{\Gg}$ can be made
into a $G$-equivariant central extension, compare (1.9.6). Suppose that such an equivariant
structure has been chosen. Then the class $c_\gamma(P)\in H^3(B, \bbC)$ can be constructed as follows. 
We have the representation
$\widetilde{\Ad}$ of $G$ on $\widetilde{\Gg}$, and therefore an extension of associated
vector bundles on $B$:
$$0\to C^\infty_B \to\widetilde{\Ad}(P) \to \Ad(P)\to 0.$$ 
Choose a connection $\nabla$ in $P$. Then we have associated linear connections
$\nabla_{\Ad}$ in $\Ad(P)$ and $\nabla_{\widetilde{\Ad}}$ in $\widetilde{\Ad}(P)$. We
also have the curvature $R_\nabla\in \Omega^2(B)\otimes\Ad(P)$. Choose a lifting $\widetilde{R_\nabla}$
of $R_\nabla$ to $\Omega^2(B)\otimes \widetilde{\Ad}(P)$, and take
$$S = \nabla_{\widetilde{\Ad}}(\widetilde{R}_\nabla)\in \Omega^3(B)\otimes\widetilde{\Ad}(P).$$
By the Bianchi identity $\nabla(R_\nabla)=0$ and so $S$ lies in the tensor product of $\Omega^3(B)$ and the subbundle
$C^\infty_B\subset \widetilde{\Ad}(P))$, i.e., it is a scalar differential form $S\in\Omega^3(B)$. 
Further, it is clear that $S$ is a closed 3-form. The class $c_\gamma(P)$ is then the class of the form $S$. 
A different choice of an equivariant structure on $\widetilde{\Gg}$ leads to change of the class of $S$
by an element from the image of $d_2$.

\vskip .2cm

(c) Let $G= GL_N(\bbC)$, so $\Gg = {\frak {gl}}_N(\bbC)$. Then $H^\bullet(\Gg)$ is the
exterior algebra on generators $\gamma_1, ..., \gamma_N$ with $\gamma_i\in H^{2i-1}(\Gg)$.
A principal $G$-bundle $P$ on $B$ is the same as a rank $N$ vector bundle $E$. In this
case each $1\otimes\gamma_i$ is transgressive, and
$c_{\gamma_i}(P)$ is the image of  $c_i(E)\in H^{2i}(B)$ under the natural
projection $H^{2i}(B)\to E_{n+1}^{0, n+1}$. Here $c_i(E)$
is the usual $i$th Chern class of $E$.

\vskip .3cm

\noindent{\bf (3.2) Other interpretations.}  Here we collect, for future use, some more
or less straightforward reformulations of the construction of $c_\gamma(P)$.  

\vskip .2cm

\noindent {\sl (a) The Chern-Weil picture.} If we choose  a connection $\nabla$ in $P$, then
the sequence (1.6.2) splits (such splitting is in fact the definition of a connection
following Atiyah). So we can identify
$$\Omega^\bullet(P)^G = \DR^\bullet(\cA_P) = \Omega^\bullet_P\otimes C^\bullet(\Ad(P)_{/B}). \leqno (3.2.1)$$
Let $R$ be the curvature of $\nabla$. Then the differential in the RHS of (3.2.1)
has the form $\partial+\nabla +i_R$, where $\partial$ is the differential in $C^\bullet(\Gg)$ and
$$i_R: \Omega^\bullet_B\otimes C^\bullet(\Gg) \to \Omega^{\bullet +2}_B\otimes
C^{\bullet -1}(\Gg)\leqno (3.2.2)$$
is the contraction with $R$. This leads to a definition of $c_\gamma(P)$ in terms of
differential forms. Namely, we have an injective and a surjective morphisms of complexes:
$$\Omega^\bullet_B = \Omega^\bullet_B\otimes C^0(\Gg)\buildrel\phi\over \hookrightarrow   
\Omega^\bullet_B\otimes C^\bullet(\Gg) \buildrel \psi\over\longrightarrow \Omega^0_B\otimes
C^\bullet(\Gg).\leqno (3.2.3)$$
Here $\psi$ is identified with the projection to $\gr^0_F$, where $F$ is the filtration from
(3.1.2). If our class $\gamma$ is basic, then it lifts uniquely to a class in
$H^n( {\operatorname{Coker}}(\phi))$, so $c_\gamma(P)$ is the image of that lifted class under
the coboundary map  corresponding to the short exact sequence
$$0\to \Omega^\bullet_B\buildrel\phi\over\hookrightarrow \Omega^\bullet_B\otimes C^\bullet(\Gg)
\to {\operatorname{Coker}}(\phi)\to 0. \leqno (3.2.4)$$

\vskip .2cm

\noindent {\sl (b) The differential graded picture.} 
Let ${\mathfrak A}$ denote the cone of the map $i : \Ad(P)\to\cA_P$
viewed as a differential graded Lie algebroid. Thus $\cA_P$ is put in degree 0, and
$\Ad(P)$ in degree $(-1)$. 
The anchor map $\alpha$ induces
the quasi-isomorphism of Lie algebroids ${\mathfrak A} @>>> {\mathcal T}_B$,
hence the map of respective universal enveloping (differential graded) algebras
${ U}({\mathfrak A}) @>>> { U}({\mathcal T}_B)=
{\mathcal D}_B$ (the latter concentrated in degree zero) which is a quasi-isomorphism. Define the map
$$
{\operatorname{DR}}^{\bullet}({\mathcal A}_P)/{\operatorname{DR}}^{\bullet}({\cT_B}) @>{\delta}>> {\operatorname{DR}}^{\bullet+1}({\mathfrak A})
\leqno (3.2.5)$$
as follows. For $X\in \Ad(P)$, denote by ${\underline X}$ the element $(X,0)$ in the cone ${\mathfrak A}$ of $i$; for $Y \in {\mathcal A}_P$, 
denote the element $(0,Y)$ simply by $Y$.
 Given a $p$-cochain $\omega$ from ${\operatorname{DR}}^{\bullet}({\mathcal A}_P$, define the cochain $\delta \omega$ by
$$
\delta \omega ({\underline X}_1, \ldots, {\underline X}_q, Y_1, \ldots, Y_r)=\omega ({\underline X}_1, Y_1, \ldots, Y_r)
\leqno (3.2.6)$$
for $q=1$ and zero for $q\neq 1$.

It is easy to see that the sequence
$$
{\operatorname{DR}}^{\bullet}({\mathcal A}_P{\operatorname{DR}}^{\bullet}({\cT_B}) @>{\delta}>>
 {\operatorname{DR}}^{\bullet+1}({\mathfrak A}) \leftarrow {\operatorname{DR}}^{\bullet+1}(\cT_B)=\Omega^{\bullet+1}_B
\leqno (3.2.7)$$
represents the boundary map
$$H^\bullet(\DR^\bullet(\cA_P)/\DR^\bullet(\cT_B))\to H^{\bullet+1}(\DR^\bullet(\cT_B)) = H^{\bullet+1}(B). \leqno
(3.2.8)$$
A basic class $\gamma$ as above defines an $n$-dimensional cohomology class $\tilde\gamma$ of
 ${\operatorname{DR}}^{\bullet}({\mathcal A}_P)/{\operatorname{DR}}^{\bullet}({\cT_B})$,
and $c_\gamma(P)$ is the image of $\tilde\gamma$ under (3.2.7).

\vskip .2cm

\noindent {\sl (c) The $\cD$-module picture.} 
Consider the short exact sequence
$$0\to C^{\geq 1}(\Ad(P)_{/B}) \to C^\bullet(\Ad(P)_{/B}) \to C^\infty_B\to 0
\leqno (3.2.9)$$
coming  from the fact that $C^\infty_B = C^0(\Ad(P)_{/B})$ is the 0th term of the relative
cochain complex. If $\mathfrak A$ is as in (b), then all three complexes in (3.2.9)
are graded $U({\mathfrak A})$-modules in the following way. Elements $Y=(0,Y), Y\in \cA$,
act via the adjoint action. Elements $\underline{X} = (X,0), X\in\Ad(P)$, act by exterior
multiplication. 
The action of $U({\mathfrak A})$ on $C^\infty_B$ is via the quasiisomorphism with $\cD_B$.

Note that (3.2.9) splits as a short exact sequence of complexes of vector bundles
but not of $U({\mathfrak A})$-modules. We will use the  corresponding connecting morphism
$$\delta: C^\infty_B\to C^{\geq 1}(\Ad(P)_{/B})[1]
\leqno (3.2.10)$$
in $D(U({\mathfrak A}))$, the derived category of differential graded $U({\mathfrak A})$-modules. 

As $\mathfrak A$ is quasiisomorphic to $\cT_B$, the DG algebra $U({\mathfrak A})$ is quasiisomorphic to
$\cD_B$, and the category $D(U({\mathfrak A}))$ is equivalent to $D(\cD_B)$. Now recall (1.5.3) that
$$H^m(B) = \Hom_{D(\cD_B)}(C^\infty_B, C^\infty_B[m]). \leqno (3.2.11)$$
On the other hand, suppose that $\Gg$ is such that $H^i(\Gg)=0$ for $0<i<n$, see (3.1.11).
Then $H^i(\Ad(P)_{/B}) = H^i(\Gg)\otimes C^\infty_B = 0$ for $0<i<n$ as well. In other words,
the complex $C^{\geq 1}(\Ad(P)_{/B})$ is acyclic in degrees $<n$ and therefore each
class  $\xi$ in its $n$th cohomology (which is isomorphic to $H^n(\Gg)\otimes C^\infty_B$) defines a
morphism in the derived category of complexes of vector bundles
$$\tilde\xi: C^{\geq 1}(\Ad(P)_{/B})\to  C^\infty_B[n]. \leqno (3.2.12)$$
Further,  ``constant'' class $\xi$, i.e., a class of the form $\gamma\otimes 1$, $\gamma\in H^n(\Gg)$, 
defines in fact a morphism in the category  $D(U({\mathfrak A}))\sim  D(\cD_B)$.
Composing $\widetilde{\gamma\otimes 1}$ with $\delta$, we get a morphism
$$C^\infty_B\to C^\infty_B[n+1], \leqno (3.2.13)$$
i.e., a class in $H^{n+1}(B)$.

\proclaim (3.2.14) Proposition. The class in $H^{n+1}(B)$ corresponding to (3.2.13)
is equal to $c_\gamma(P)$. 

\noindent {\sl Proof:} This follows directly from the definitions (in fact, we could take
(3.2.13) as the definition of $c_\gamma(P)$). Indeed, the morphism in the derived category
from the cohomology of a quotient complex such as $C^\infty_B$ to the homology of
a  subcomplex such as $C^{\geq 1}(\Ad(P)_{/B})$ acyclic up to degree $n$, is precisely the differential
$d_{n+1}$ in the corresponding spectral sequence. \qed

\vskip .3cm

\noindent {\bf (3.3) Infinite-dimensional groups.} Slightly reformulating the approach of
K.-T. Chen [C], we introduce the following definition.

\proclaim (3.3.1) Definition. A differentiable space is an ind-object in the category
of $C^\infty$-manifolds.

For background on ind-objects, see [D2]. Thus a differentiable space $M$ is a formal limit
$``\lim\limits_{\longrightarrow}{}''{}_{\alpha\in A} M_\alpha$ of (finite-dimensional)
$C^\infty$-manifolds. In particular, $M$ defines a functor
$$S\mapsto M(S) = C^\infty(S, M) = \lim_{\longrightarrow} C^\infty(S, M_\alpha)
\leqno (3.3.2)$$
on such manifolds and can in fact be identified with this functor. In practice, however,
we will identify
$M$  with the set $M(pt) = \lim\limits_{\longrightarrow} M_\alpha$
with (3.2.2) providing an additional structure on this set (description of what it means
for an element of this set to vary in a smooth family). \

For a differential space $M$ we define (compare [C]) the space of $p$-forms (in particular, of
$C^\infty$-functions) on $M$ by
$$\Omega^p(M) = \lim_{\longleftarrow} \Omega^p(M_\alpha).\leqno (3.3.3)$$
For a point $m\in M(pt)$ the tangent space $T_mM$ is defined  by
$$T_mM = \lim_{\longrightarrow} \, T_sS,\leqno (3.3.4)$$
where the limit is taken over $C^\infty$-maps $(S,s)\to (M,m)$.

A differentiable group $G$ is a group object in the category of differentiable spaces.
For such a group the space $\Gg= T_eG$ is a Lie algebra in the standard way. 

\vskip .2cm

\noindent {\bf (3.3.5) Examples.} {\sl (a) Groups of diffeomorphisms.} Let $\Sigma_0$ be
a compact oriented $C^\infty$-manifold of dimension $d$. Then we have a differentiable group
$G = \operatorname{Diffeo}(\Sigma_0)$ of orientation preserving diffeomorphisms.
The corresponding functor (3.3.2) is as follows. A smooth map $S\to \operatorname{Diffeo}(\Sigma_0)$
is a diffeomorphism of $S\times\Sigma_0$ preserving the projection to $S$. The Lie algebra
of this group is $\operatorname{Vect}(\Sigma_0)$, the algebra of $C^\infty$ vector fields.

\vskip .1cm

\noindent{\sl (b) Gauge groups.} Let $\Sigma_0$ be as before and $E_0$ be a $C^\infty$
complex vector bundle on $\Sigma_0$ Then we have the differentiable group $\Aut(E_0)$
of $C^\infty$-automorphisms of $E_0$ (the differentiable structure defined similarly to (a)).
Its Lie algebra is $\End(E_0)$. 

\vskip .1cm

\noindent {\sl (c) Atiyah groups.} Let $\Sigma_0, E_0$ be as before. The Atiyah group $AT(\Sigma_0, E_0)$
consists of pairs $(\phi,f)$, where $\phi$ is an orientation preserving diffeomorphism
of $\Sigma_0$, and $f: \phi^*E_0\to E_0$ is an isomorphism of vector bundles. Thus we have an
extension of differentiable groups:
$$1\to\Aut(E_0) \to AT(\Sigma_0, E_0) \to \operatorname{Diffeo}(\Sigma_0)\to 1.$$
The Lie algebra of $AT(\Sigma_0, E_0)$ is $\cA_{E_0}(\Sigma_0)$, the algebra of
global $C^\infty$-sections of the Atiyah Lie algebroid. 

In fact, the Atiyah group is a particular case of constructions in [Mac], \S 1.4.4-7, namely
the group of bisections of the frame groupoid of $E_0$. 

More generally, one can replace the vector bundle in Examples (b), (c) by a principal
bundle with an arbitrary structure Lie group. In this paper we will
be interested in the vector bundle case and will concentrate on the Example (c) as
the most general. 

\vskip .2cm

Let us now describe a class of principal bundles with structure groups as in (c).
Suppose $q: \Sigma\to B$ is a smooth fibration with compact oriented fibers of dimension $d$.
Suppose that $B$ is connected. Then all the fibers $\Sigma_b = q^{-1}(b), b\in B$, are diffeomorphic
to each other. Let $\Sigma_0$ be one such fiber. Futher, let $E$ be a smooth $\bbC$-vector bundle
on $\Sigma$ and $E_b= E|_{\Sigma_b}$. Then for different $b$ the pairs $(\Sigma_b, E_b)$ are
isomorphic, in particular, isomorphic to $(\Sigma_0, E_0)$. Let $G= AT(\Sigma_0, E_0)$.
We have then the principal $G$-bundle
$$\rho: P=P(\Sigma/B, E)\to B
\leqno (3.4.5)$$
whose fiber $P_b=\rho^{-1}(b), b\in B$, consists of isomoprhisms of pairs 
$(\Sigma_0, E_0)\to (\Sigma_b, E_b)$.

\vskip .1cm

For any differentiable $G$-bundle $P$ over a finite-dimensional base $B$ the
Atiyah algebra $\cA_P$ can be defined by (1.6.3). In the example where $G=AT(\Sigma_0, E_0)$
and $P= P(\Sigma/B, E)$, this gives
$$\cA_{P(\Sigma/B, E)} = q_* \cA_E \leqno (3.4.6)$$
(the sheaf-theoretic direct image of the Atiyah algebra of $E$). 

\vskip .3cm

\noindent {\bf (3.5) The first Chern class.} Let $q: \Sigma\to B$ and $E$ be as before,
so that we have a principal bundle $P=P(\Sigma/B, E)\to B$ with structure group
$G=AT(\Sigma_0, E_0)$. AS the corresponding Lie algebra $\Gg= \cA_{E_0}(\Sigma_0)$
consists of global sections of the Atiyah Lie algebroid of $\Sigma_0$, we have the embeddings
$$\Gg \hookrightarrow \cD(E_0) \hookrightarrow  {\mathfrak{gl}} (\cD(E_0)).\leqno (3.5.1)$$
By Corollary 2.3.4, ${\mathfrak gl} (\cD(E_0))$ has a unique continuous (in the Frechet topology)
cohomology class $c$ in degree $d+1$. We denote by $\gamma$ the restriction of $c$ to $\Gg$.

\proclaim (3.5.2) Proposition-Definition. The class $1\otimes\gamma$ is transgressive,
so $d_{d+2}(1\otimes\gamma)$ is defined. 
 Moreover,
the Lie algebroid $\cA_P$ naturally embeds into a bigger Lie algebroid
$$0\to q_* ({\mathfrak gl}(\cD_{\Sigma/B}(E)))\to \cA_{\Sigma/B, E} \buildrel \alpha\over
\rightarrow \cT_B\to 0. \leqno (3.5.3)$$
Here the fibers of $\Ker(\alpha)$ are Lie algebras isomorphic to ${ \mathfrak {gl}}(\cD_{\Sigma_0, E_0}))$
 via an isomoprhism
defined uniquely up to an inner automorphism and thus satisfy the sphericity condition (3.1.5) with $n=d+1$.
Therefore $d_{d+2}(1\otimes\gamma)$ admits a canonical lifting to a class
 $$C_1(q_*E) := c_\gamma(P) \in H^{d+2}(B, \bbC).$$

This class will be the main object of study in the rest of the paper. 

\vskip .1cm

To prove Proposition 3.5.2, 
it is enough to describe the construction of $\cA_{\Sigma/B, E}$, as the Hochschild-Serre
spectral sequence for it maps into the analogous sequence for $\cA_P$.

\vskip .2cm

\noindent {\bf (3.5.4) Construction of $\cA_{\Sigma/B, E}$.}
We start with the Atiyah Lie algebroid on
$\Sigma$:
$$0\to \End(E)\buildrel i\over\rightarrow\cA_E\buildrel\alpha\over
\rightarrow \cT_\Sigma\to 0.\leqno (3.5.5)$$

Let ${ U}({\mathcal A})_{E/B}$
denote the centralizer of $q^{-1} C^\infty_B$ in
${ U}({ A_E})$.
Let $F_\bullet{ U}({\mathcal A}_E)$ denote the filtration defined by
$F_{-1}{ U}({\mathcal A}_E) = 0$, $F_0{ U}({\mathcal A}_E) =
{ U}({\mathcal A})_{E/B}$ and
$[F_i{ U}({\mathcal A}_E), q^{-1} C^\infty_B]\subseteq
F_{i-1}{ U}({\mathcal A}_E)$. Then, $F_1{ U}({\mathcal A}_E)$
is a Lie algebra under the commutator, ${ U}({\mathcal A})_{E/B}$
is a Lie ideal in $F_1{ U}({\mathcal A}_E)$, and there is an exact
sequence
$$
0 \to { U}({\mathcal A})_{E/B} \to
F_1{ U}({\mathcal A}_E) \to  q^{-1}{\mathcal T}_B \to  0
\leqno (3.5.6)$$
exhibiting $F_1{ U}({\mathcal A}_E)$ as a transitive
$q^{-1} C^\infty_B$-algebroid.

 The inclusion ${ \cA}_{E} \to
{\mathcal D}_{\Sigma} (E)$ induces the surjective map
${ U}({\mathcal A}_{ E})\to {\mathcal D}_{\Sigma} (E)$
with kernel being the ideal generated by the relation which identifies
$1\in C^\infty_\Sigma \subset{ U}({\mathcal A}_{E})$ with
$1\in\shEnd_{C^\infty_\Sigma}({ E})\subset{\mathcal A}_{ E}$.
The exact sequence (3.5.6)  reduces to the exact sequence
$$0 \to {\mathcal D}_{\Sigma/B,{ E}}\to
F_1{\mathcal D}_{\Sigma,{ E}} \to q^{-1}{\mathcal T}_B \to  0.
\leqno (3.5.7)$$

Replacing ${ E}$ by its tensor product by the trivial bundle of rank $r$ in the above
example, (3.5.7)  can be rewritten as
$$
0 \to  {\mathfrak gl}_r({\mathcal D}_{\Sigma/B, { E}}) \to
F_1{\mathfrak gl}_r({\mathcal D}_{\Sigma, { E}}) \to  q^{-1}{\mathcal T}_B \to 0.
\leqno (3.5.8)$$
Taking the limit over inclusions ${\mathfrak gl}_r \to {\mathfrak gl}_{r+1}$
we obtain a $q^{-1}C^\infty_B$-algebroid
$$0\to  {\mathfrak gl}({\mathcal D}_{\Sigma/B, { E}}) \to
{\mathcal A}_{q, { E}} \to  q^{-1}{\mathcal T}_B \to  0. \leqno (3.5.9)$$
Put 
$${\mathfrak A}_{q, {E}}={\operatorname{Cone}}\biggl\{{\mathfrak gl}({\mathcal D}_{\Sigma/B, { E}}) \to
{\mathcal A}_{q, { E}}\biggr\}. \leqno (3.5.10)$$
There are quasi-isomorphisms 
$$
{\mathfrak A}_{q, {E}} \to q^{-1}C^\infty_B, \quad 
{ U}_{q^{-1} C^\infty_B}({\mathfrak A}_{q, { E}})\to q^{-1}{\mathcal D}_B. $$

Taking the direct image of (3.5.9)  under $q$ and pulling back by
the canonical map ${\mathcal T}_B \to q_* q^{-1}{\mathcal T}_B$ we obtain
the transitive
(because ${ R}^1\pi_*{\mathfrak gl}_r({\mathcal D}_{\Sigma/B, { E}}) = 0$)
Lie algebroid on $B$:
$$
0 \to {\mathcal G} \to
{\mathcal A}_{\Sigma/B,  {E}} \to {\mathcal T}_B \to 0, \leqno (3.5.11)$$
where ${\mathcal G} = q_*{\mathfrak gl}({\mathcal D}_{\Sigma/B, { E}})$, as we wanted.
Let ${\mathfrak A}_{\Sigma/B, { E}}$ denote the differential graded Lie algebroid
on $B$ equal to the cone of the inclusion
${\mathcal G} \to {\mathcal A}_{\Sigma/B,  E}$. Let
$$\delta_{\Sigma/B} : C^\infty_B \to  C_+({\mathcal G})[1]
\leqno (3.5.12)$$
denote the connecting homomorphism; it is a morphism in the derived category of differential graded modules over the universal
enveloping (differential graded) algebra ${ U}({\mathfrak A}_{\Sigma/B})$.

\vskip .3cm

\noindent {\bf (3.6) Smooth cohomology and characteristic classes.}
A more traditional way of getting characteristic classes of principal $G$-bundle is by using group
cohomology classes of $G$. Let us present a framework which we will then compare with
the Lie algebra framework above.

Let $S$ be a topological space and $\cF$ be a 
sheaf of abelian groups on $S$. We denote by $\Phi^\bullet(\cF)$ the standard Godement resolution of
$\cF$ by flabby sheaves. Thus $\Phi^0(\cF) = DS(\cF)$ is the sheaf of (possibly discontinuous)
sections of the (etale space associated to) $\cF$, and $\Phi^{n+1}(\cF) = DS(\Phi^n(\cF))$. 
In this and the next sections we write $R\Gamma(S, \cF)$ for the complex of global sections
$\Gamma(S, \Phi^\bullet(\cF))$.

Let $G$ be a differentiable group and $B_\bullet G$ be its classifying space. Thus
$B_\bullet G= (B_nG)_{n\geq 0}$ is a simplicial object in the category of differentiable spaces
with $B_nG=G^n$, and the face and degeneracy maps given by the standard formulas. We define the smooth
cohomology of $G$ with coefficients in $\bbC^*$ to be
$$H^n_{sm}(G, \bbC^*) = \bbH^n(B_\bullet G, C^{\infty *}). \leqno (3.6.1)$$
Here the hypercohomology on the right is defined as the cohomology of the double complex 
whose rows  are  the
complexes $R\Gamma(B_nG, C^{\infty *}_{B_nG})$ and the differential between the neighboring slices
coming from the simplicial structure on $B_\bullet G$. 
This is a version of the Segal cohomology theory for topological groups ([F], p. 305).
In particular, we have a spectral sequence
$$H^i(B_n G, \bbC^*) \Rightarrow H^{i+n}_{sm}(G, \bbC^*).\leqno (3.6.2))$$
We will use some other natural (complexes of) sheaves on $B_\bullet G$ to get natural
cohomology theories for $G$. For example, the   Deligne cohomology
$$H^n_{sm}(G, \bbZ_D(p)) = \bbH^n(B_\bullet G, \bbZ_D(p)), \leqno (3.6.3)$$
where for any differentiable space $M$ we set
$$\bbZ_D(p) = \biggl\{ \underline{\bbZ}_M\to \Omega^0_M\to\Omega^1_m\to ... \to \Omega^{p-1}_M\biggr\},
\leqno (3.6.4)$$
with $\underline{\bbZ}_M$ put in degree 0, compare [Bry].

\vskip .2cm

Let $B$ be a $C^\infty$-manifold and ${{\cU}} = \{U_i\}_{i\in I}$ be an open covering of $B$.
We denote by $N_\bullet{{\cU}}$ the simplicial nerve of ${{\cU}}$, i.e., the
simplicial manifold with
$$N_n{{\cU}} = \coprod_{i_0, ..., i_n} U_{i_0} \cap ...\cap U_{i_n}.\leqno (3.6.5)$$
As well known, $N_\bullet{{\cU}}$ is homotopy equivalent to $B$, so for any sheaf $\cF$
on $B$ we have
$$\bbH^i(N_\bullet {\cU}, \cF_\bullet) = H^i(B, \cF), \leqno (3.6.6)$$
where $\cF_\bullet$ is the natural sheaf on $N_\bullet{{\cU}}$ whose
$n$th component is the sheaf on (3.6.6) formed by the restrictions of $\cF$. 

Let $\rho: P\to B$ be a principal $G$-bundle and suppose that $P$ is trivial on each $U_i$.
Then a collection of trivializations (i.e., sections) $\tau = (\tau_i: U_i\to P)$
gives a morphism of simplicial differentiable spaces 
$$u_\tau: N_\bullet{\cU}\to B_\bullet G.\leqno (3.6.7)$$
Given a class $\beta\in H^n_{sm}(G, \bbC^*)$, we define the characteristic class
$${\frak {c}}_\beta(P) = u_\phi^* (\beta)\in H^n(B, C^{\infty *}_B). \leqno (3.6.8)$$
Similarly one can define characteristic classes  corresponding
to group cohomology classes with values in the Deligne cohomology.


\vskip .3cm

\noindent {\bf (3.7) Integrality and integrability.} Let $G$ be as in (3.6), and
$\Gg$ be the Lie algebra of $G$. We  construct the ``derivative" map
$$\partial: H^n_{sm}(G, \bbC^*) \to H^n_{Lie}(\Gg, \bbC).\leqno (3.7.1)$$
To do this, we first remark that for any topological space $S$, any sheaf of abelain groups $
\cF$
on $S$ and any point $s_0\in S$ we have a natural morphism of complexes
$$\epsilon_{s_0}: R\Gamma(S, \cF) \to \cF_{s_0}, $$
where $\cF_{s_0}$ is the stalk of $\cF$ at $s_0$. To construct $\epsilon_{s_0}$, we first project
$R\Gamma(S, \cF) = \Gamma(S, \Phi^\bullet(\cF))$ to its 0th term $\Gamma(S, \Phi^0(\cF))$ which, by
definition, is the space of all sections $\phi = (s\mapsto \phi_s)$ of the etale space of $\cF$. Thus
any such $\phi$ is a rule which to any point $s\in S$ associates an element of $\cF_s$. We define $\epsilon_{s_0}$ by further
mapping any $\phi$ as above to $\phi_{s_0}\in\cF_{s_0}$. 

We now specialize to $S= B_mG=G^m$, to $s_0 = e_m :=(1, ..., 1)$ and to $\cF= C^{\infty *}_{S}$. We get a morphism
from the double complex 
$$\{R\Gamma(B_mG, C^{\infty *}_{B_mG})\}_{m\geq 0}\leqno (3.7.2)$$
to the complex of stalks
$$\bbC^* \to C^{\infty *}_{G, e_1} \to C^{\infty *}_{G\times G, e_2} \to ...\leqno (3.7.3)$$
An $n$-cocycle in (3.7.2) gives thus a germ of a smooth function 
$$\xi = \xi(g_1, ..., g_n): G^n\to \bbC^*$$
satisfying the group cocycle equation (on a neighborhood of $e_{n+1}$ in
$G^{n+1}$). Similarly to [F], p. 293, one associates
to $\xi$ a Lie algebra cocycle $\partial(\xi)\in C^{n}(\Gg)$ by
$$\partial(\xi)(x_1, ..., x_n) = {d\over dt} \log \xi(\exp(t x_1), ..., \exp (t x_n))\biggr|_{t=0}. \leqno (3.7.4)$$

\vskip .2cm
 
A Lie algebra
cohomology class $\gamma\in H^n(\Gg, \bbC)$ will be called integrable, if it lies in the image of the map
$\partial$ from (3.7.1).  Consider the exponential exact sequence (1.10.1) 
of sheaves on $B$:
 and ts coboundary map
 $\delta_n$ from (1.10.2).  The intuition with 
determinantal $d$-gerbes (0.4) suggests the following.

\proclaim (3.7.5) Conjecture. (a) The class $\gamma\in H^{d+1}(\cA_{E_0}(\Sigma_0))$ constructed in (3.5)
is integrable and comes from a natural class $\beta\in H^{d+1}_{sm}(AT(\Sigma_0, E_0), \bbC^*)$
(the ``higher determinantal class").
\hfill\break
(b) Further, for any $q: \Sigma\to B$ and $E$ as above, the class $C_1(q_*E) = c_\gamma(P) \in H^{d+2}(B, \bbC)$ is
integral and is the image of the following class in the integral cohomology:
$$\delta_{d+1}({\frak{c}}_\beta(P))\in H^{d+2}(B, \bbZ).$$

This conjecture holds for $d=1$ (i.e., for the case of a circle fibration). We will verify this
in Section 5. In general,  the property (b) seems to follow from (a) in virtue of some compatibility
result between group cohomology classes with coefficients in $\bbC^*$ and Lie algebra
cohomology classes with coefficients in $\bbC$. Here we present a $d=1$ version of such a result.

 Let $G$ be a
differentiable group with Lie algebra $\Gg$. Let $\beta\in H^2_{sm}(G, \bbC^*)$
and $\gamma = \partial(\beta)\in H^2_{\Lie}(\Gg, \bbC)$
be the derivative of $\beta$. Suppose $\beta$ 
is represented by an extension of differentiable groups
$$1\to \bbC^* \to \widetilde{G}\to G\to 1,\leqno (3.7.6) $$
whose Lie algebra is the extension
(3.1.9) representing $\gamma$.
Let $\rho: P\to B$ be a principal $G$-bundle over
a $C^\infty$-manifold $B$. Then we have the characteristic class
$c_\gamma(P)\in H^3(B, \bbC)$ (lifting to $H^3$ well defined because $\widetilde{\Gg}$
is a ${G}$-module via the adjoint representation of $\widetilde{G}$, see Example 3.1.8(b)). 
On the other hand, $\beta$ gives rise to a class ${\frak{c}}_\beta(P)\in H^2(B, C_B^{\infty *})$, 
see (3.6.8).

 \proclaim (3.7.7) Proposition. In the above situation
$c_\gamma(P)\in H^3(B, \bbC)$ is the image of $\delta_2({\frak {c}}_\beta(P))\in H^3(B, \bbZ)$
under the natural homomorphism from the integral to the complex cohomology.

\noindent {\sl Proof:} This follows from the result of Brylinski (1.10.4) using
Example 3.1.8(b) and an obvious generalization of (1.10.3) to differentiable groups. 
 \qed

\vskip .2cm

We further conjecture the existence of the natural ``deloopings" of the higher Chern classes as well, i.e.,
the existence of  classes
$$\beta_m \in H^{d+2m}_{sm}(AT(\Sigma_0, E_0), \bbZ_D(m)), \quad m\geq 1, \leqno (3.7.8)$$
which then give characteristic classes in families:
$$C_m(q_*E) \in H^{d+2m}(B, \bbZ_D(m)).\leqno (3.7.9)$$


\vfill\eject

\centerline {\bf 4. The Real Riemann-Roch.}

\vskip 1cm

Here is the main result of the present paper.

\proclaim (4.1) Theorem. Let $q:\Sigma\to B$ be a $C^\infty$ fibration  with compact oriented fibers of
dimension $d$. Let $E$ be a complex $C^\infty$ vector bundle on $\Sigma$. Then:
$$C_1(q_*E) = \int_{\Sigma/B}\biggl[ ch(E)\cdot \Td({\mathcal T}_{\Sigma/B})\biggr]_{2d+2}
\quad\in\quad 
H^{d+2}(B, {\mathbb C}).)$$

The proof consists of several steps.

\vskip .3cm

\noindent {\bf (4.2) A $\cD$-module interpretation of $C_1$ using $\cA_{\Sigma/B, E}$.}
We use the notations of (3.5.3-12) and introduce the following abbreviations:
$$\cG = q_*({\mathfrak gl} (\cD_{\Sigma/B}(E))).\leqno (4.2.1)$$
This is a bundle of infinite-dimensional Lie algebras on $B$.
$${\mathfrak A} = {\mathfrak A}_{q,E}.\leqno (4.2.2)$$
This is a DG Lie algebroid on $\Sigma$ quasiisomorphic to $\cT_\Sigma$. 
$$U{\mathfrak A} = U_{q^{-1} C^\infty_B} ({\mathfrak A}_{q,E}). \leqno (4.2.3)$$
This is a sheaf of DG-algebras on $\Sigma$ quasiisomorphic to $q^{-1}\cD_B$. 

Now, $U{\mathfrak A}$ acts on $C_+({\mathfrak gl}(\cD_{\Sigma/B, E}))_B$, the reduced
relative Lie cochain complex. Further, it acts on the relative Hochschild and cyclic
complexes. In the same spirit as in (3.2)(c), elements $Y=(0,Y), Y\in\cA_{q,E}$, act
via the adjoint action. Elements of the form  $\underline{X}=(X,0)$ act via the 
``insertion operators''
$$\iota_X(a_0\otimes ...\otimes a_p) = \sum_{i=0}^p (-1)^i
a_0\otimes ... \otimes a_i\otimes X\otimes a_{i+1}\otimes ... \otimes a_p.\leqno (4.2.4)$$
Denoting by $b,B$ the standard  operators on Hoschshild cochains, see [L],
we have 
$$[b, \iota _X] = \ad(X), \quad [b, \iota_X]=0.$$
Therefore $U{\GA}$ acts on both the Hochschild and the cyclic complexes. This
action extends to the completions described in (2.4). Further, the morphisms
$\mu_\cD, \nu_\cD$ from (2.4.4-5) are in fact morphisms in $D(U\GA)$. Indeed,
there is a spectral sequence
$$E_2^{pq} = \Ext^p_{q^{-1}\cD_B}\bigl( \underline H^q(\widehat{\Hoch}^\bullet(
\cD_{\Sigma/B}(E))), C^\infty_B\bigr) \Rightarrow
\Ext^{p+q}_{U\GA} (\bigr( \widehat{\Hoch}^\bullet(\cD_{\Sigma/B}(E)), C^\infty_B),
\leqno (4.2.5)$$
and similarly for the cyclic complex.
The map $\mu_\cD$ defines an element of $E_2^{0d}$, and $E_2^{pq}=0$ for $q<d$, so
$\mu_\cD$ gives rise to a well defined class in $\Ext^d$ on the RHS of (4.2.5).
Similarly for $\nu_\cD$. 

Let 
$$\alpha: C_+({\mathfrak gl}(\cD_{\Sigma/B}(E)))_B\to q^{-1} C^\infty_B[2d]
\leqno (4.2.6)$$
denote the composition
$$C_+({\mathfrak gl}(\cD_{\Sigma/B}(E)))_B\buildrel\beta\over\longrightarrow
CC_\bullet(\cD_{\Sigma/B}(E))_B[1] \to \widehat{CC}_\bullet(\cD_{\Sigma/B}(E))_B[1]
\buildrel \nu_\cD[1]\over\longrightarrow q^{-1} C^\infty_B[2d+1].\leqno (4.2.7)$$
Here the first morphism $\beta$ is the standard map from the Lie algebra chain
complex to the cyclic complex, see [L], (10.2.3). 

It is checked directly that $\beta$ commutes with the operators $\iota_X$, so it is
$U\GA$-invariant. Therefore, all maps in (4.2.7) and the map $\alpha$ are
morphisms in $D(U\GA)$. 

Let us now thake the $C^\infty$ direct image and define the morphism
$$\int_{\Sigma/B}\alpha: C_+(\cG)_B\to C^\infty_B[d] \leqno (4.2.8)$$
as the composition
$$C_+(\cG)_B \to q_* C_+{\mathfrak gl}(\cD_{\Sigma/B}(E))_B \buildrel\sim\over
\longrightarrow Rq_*C_+ {\mathfrak gl} (\cD_{\Sigma/B}(E))_B \buildrel\alpha\over\longrightarrow$$
$$\to Rq_*q^{-1} C^\infty_B[2d+1]\buildrel \int_{\Sigma/B}\over\longrightarrow C^\infty_B[d+1].
\leqno (4.2.9)$$
Here the last map is the integration over the relative (topological)
fundamental class of $\Sigma/B$. 

Consider the composition
$$C^\infty_B\buildrel\delta_{\Sigma/B}\over\longrightarrow
C_+(\cG)_B[1]\buildrel \int_{\Sigma/B}\alpha\over\longrightarrow C^\infty_B[d+2],
\leqno (4.2.10)$$
where $\delta_{\Sigma/B}$ is as in (3.5.12). As both maps in (4.2.10) are morphisms in
$D(\cD_B)$, the composition (denote it $C$) is an element
$$C\in \Ext^{d+2}_{\cD_B}(C^\infty_B, C^\infty_B) = H^{d+2}(B, \bbC). $$

\proclaim (4.2.11) Proposition. We have $C= C_1(q_*E)$. 

\noindent {\sl Proof:} This follows from the interpretation of
$C_1(q_*E) = c_\gamma(P(\Sigma/B, E))$ given in (3.2) (b), (c) and from
the compatibility of the Atiyah algebroid of $P(\Sigma/B)$ with
$\cA_{\Sigma/B, E}$. 

\vskip .3cm

\noindent {\bf (4.3) A local RRR in the total space.}
Proposition 4.2.11 reduces the RRR to the following ``local" statement taking place
in the total space $\Sigma$. 

\proclaim (4.3.1) Theorem. Let $\xi$ be the morphism in $D(q^{-1} \cD_B)$ defined as the
composition
$$q^{-1}\cO_B\to C_+{\mathfrak gl} (\cD_{\Sigma/B, E}))[1] \to q^{-1}\cO_B[2d+2].$$
Then the class in
$$\Ext^{2d+2}_{q^{-1}\cD_B}(q^{-1} \cO_B, q^{-1} \cO_B) = H^{2d+2}(\Sigma, \bbC),$$
corresponding to $\xi$, is equal to
$$\biggl[ch(E)\cdot \Td(\cT_{\Sigma/B})\biggr]_{2d+2}.$$

We now concentrate on the proof of Theorem 4.3.1.  First, we remind the definition of periodic cyclic
homology [L]. Let $A$ be an associative algebra. The ``negative" cyclic complex of $A$ is defined,
similarly to (2.1.4), as
$$CC_\bullet^-(A) = \operatorname{Tot} \biggl\{ \Hoch_\bullet(A) \buildrel N\over\longrightarrow
\Hoch_\bullet(A)\buildrel 1-\tau\over\longrightarrow \Hoch_\bullet(A) \to ...\biggr\}
\leqno (4.3.2)$$
Here, the grading of the copies of $\Hoch_\bullet(A)$ in the horizontal direction
goes in increasing integers 0,1,2 etc. So $CC_\bullet^-(A)$ 
s a module over the formal Taylor series ring $\bbC [[u]]$ where
$u$ has degree $(-2)$. The original cyclic complex is a module over  the polynomial
ring $\bbC[u^{-1}]$. Finally, the periodic cyclic complex $CC_\bullet^{per}(A)$ is
obtained by merging together $CC_\bullet(A)$ and $CC_\bullet^{-}(A)$ into
one double complex which is repeated 2-periodically both in the positive and negative 
horizontal directions.
In other words,
$$CC_\bullet^{per}(A) = CC_\bullet^-(A)\otimes _{\bbC[[u]]} \bbC ((u)).\leqno (4.3.3)$$
We extend these construction to other situations (see \S 2) where the tensor
products are understood in the sense of various completions. In particular, the morphism
$\nu_D$ from (2.4.5) extends to morphisms
$$\nu_D^-: CC^-_\bullet(\cD_{\Sigma/B, E})\to q^{-1}C^\infty_B[2d] [[u]], \quad 
\nu_D^{per}: CC_\bullet^{per}(\cD_{\Sigma/B, E}) \to  q^{-1}C^\infty_B[2d]((u)). $$
These morphisms include into the commutative diagram
$$
\begin{CD}
CC^-_\bullet ({\mathcal D}_{\Sigma/B , { E}}) @>>> CC^{per}_\bullet ({\mathcal D}_{\Sigma/B, { E}})
@>>> CC({\mathcal D}_{X/S, {\mathcal E}})[2] \\
@V{\nu_{\mathcal D}^-}VV @V{\nu_{\mathcal D}^{per}}VV
@VV{\nu_{\mathcal D}}V \\
{\mathcal O}_S[2d][[u]] @>>> {\mathcal O}_S[2d]((u))
@>{Res_{u=0}}>> {\mathcal O}_S[2d+2]
\end{CD} \leqno (4.3.4)
$$
We now want to reduce Theorem 4.3.1  to  the following statement.

\proclaim (4.3.5) Theorem. The composition
$$C^\infty_B \buildrel +1 \over\rightarrow CC_\bullet^{per}(\cD_{\Sigma/B, E}) \buildrel
\nu_\cD^{per}\over\longrightarrow C^\infty_B[2d]((u))$$
is equal to
$$\sum_{i=0}^\infty u^i \cdot \bigl[\ch(E)\Td(\cT_{\Sigma/B})\bigr]_{2(d-i)}.$$

Indeed, suppose we know Theorem 4.3.5. To prove Theorem 4.3.1, it would then be sufficient to
prove that the composition
$$q^{-1} C^\infty_B\to C_+({\frak{gl}}(\cD_{\Sigma/B}))[1] \to CC_\bullet (\cD_{\Sigma/B, E})[2]
\leqno (4.3.6)$$
is equal to the composition
$$q^{-1}C^\infty_B  \buildrel +1\over\rightarrow CC_\bullet^{per} (\cD _{\Sigma/B, E}) \to
CC_\bullet (\cD_{\Sigma/B, E})[2],\leqno (4.3.7)$$
as the latter one is related to Chern and Todd via (4.3.5). 
In order to perform the comparison, let $K$ be the cone of the inclusion
$C_+({\frak{gl}}(\cD_{\Sigma/B, E})\to C_\bullet({\frak{gl}}(\cD_{\Sigma/B, E}))$, so that we have 
a quasi-isomorphism $K\to q^{-1}C^\infty_B$ as well as an isomorphism of distinguished triangles
$$
\begin{CD}
C_\bullet({\mathfrak gl}({\mathcal D}_{\Sigma/B, { E}})) @>>> K @>>> C_+({\mathfrak
gl}({\mathcal D}_{\Sigma/B, { E}}))[1] \\
@VVV @VVV @VVV \\
C({\mathfrak gl}({\mathcal D}_{\Sigma/B, {\mathcal E}})) @>>> q^{-1}C^\infty_B
@>>> C_+({\mathfrak gl}({\mathcal D}_{\Sigma/B,  E}))[1]
\end{CD}
$$
(with the top row a short exact sequence of complexes).
Notice now that there is a morphism of distinguished triangles
$$
\begin{CD}
C_\bullet ({\mathfrak gl}({\mathcal D}_{\Sigma/B, {E}})) @>>> K @>>> C_+({\mathfrak
gl}({\mathcal D}_{\Sigma/B, {E}}))[1] \\
@VVV @VVV @VVV \\
CC_\bullet^-({\mathcal D}_{\Sigma/B, { E}}) @>>> CC_\bullet^{per}({\mathcal D}_{\Sigma/B, { E}}) @>>>
CC_\bullet ({\mathcal D}_{\Sigma/B, { E}})[2]
\end{CD}
$$
It remains to notice further that the diagram
$$
q^{-1}C^\infty_B@<<< K @>>> CC_\bullet ^{per}({\mathcal D}_{\Sigma/B, { E}})
$$
represents the morphism $C^\infty_B @>{1}>> CC^{per}_\bullet ({\mathcal
D}_{\Sigma/B, { E}})$ in the derived category and the proof is finished.

\vskip .3cm

\noindent {\bf (4.4) Proof of Theorem 4.3.5.} This statement can be deduced from the results of
[NT] on the cohomology of the Lie algebras of formal vector fields and formal matrix functions. 
We recall the setting of [NT] which extends that of the Chern-Weil definition of characteristic
classes. Recall that the latter provides a map
$$S^\bullet[[{\mathfrak{h}}_0]]^{H_0} \to H^{2\bullet}(\Sigma, \bbC),\leqno (4.4.1)$$
where $H_0 = GL_d(\bbC)\times GL_r(\bbC)$ with $r= \rk(E)$, while
$\mathfrak{h}_0$ is the Lie algebra of $H_0$, i.e.,
${\mathfrak{gl}}_d(\bbC) \oplus {\mathfrak{gl}}_r(\bbC)$.  To be precise, the elementary
symmetric functions of the two copies of $\mathfrak{gl}$ are mapped to the Chern
classes of $\cT_{\Sigma/B}$ and $E$.  

In [NT], this construction was generalized in the following way. Let $k=\dim(B)$, and
  ${\widehat {\mathfrak g}}$ be the Lie algebra of formal differential operators of the form 
$$\sum_{i=1}^k P_i(y_1, \ldots, y_k){\frac{\partial}{\partial{y_i}}}+\sum_{j=1}^d Q_j(x_1, \ldots, x_d, y_1, \ldots, y_k){\frac{\partial}{\partial{x_i}}}+
R(x_1, \ldots, x_d, y_1, \ldots, y_k)\leqno (4.4.2)$$
where $P_i$, $Q_j$ are formal power series, $R(x)$ is an $r\times r$ matrix whose entries are power series.
Thus ${\widehat {\mathfrak g}}$ is the formal version of the relative Atiyah algebra. 
Consider the Lie subalgebra $\mathfrak{h}$ of fields such that all $P_i$ and $Q_j$ are of degree one and all entries of $R$ are of degree zero. 
We can identify this subalgebra with 
$$
{\mathfrak h}={\mathfrak {gl}}_d ({\mathbb C}) \oplus {\mathfrak {gl}}_k ({\mathbb C})\oplus {\mathfrak {gl}}_r({\mathbb C})
$$
Let 
$$
{H} ={\operatorname{GL}}_d ({\mathbb C}) \times {\operatorname{GL}}_k ({\mathbb C})\times {\operatorname{GL}}_r ({\mathbb C})
$$
be the corresponding Lie group. Thus $({\widehat{\mathfrak g}}, H)$ form a Harish-Chandra pair.
Following the ideas of  ``formal geometry" (or ``localization") of Gelfand and Kazhdan, one sees that
every $({\widehat{\mathfrak g}}, H)$-module $L$ induces a sheaf ${\mathcal L}$ on $\Sigma$.
Similarly, a complex $L^\bullet$ of modules gives rise to a complex of sheaves $\cL^\bullet$. A complex $L^\bullet$ of modules is
called   {\em homotopy constant}, i.e. the action of ${\widehat {\mathfrak g}}$ extends to an action of the differential graded Lie algebra 
$({\widehat {\mathfrak g}}[\epsilon], {\frac{\partial}{\partial{\epsilon}}})$. Here $\epsilon$ is a formal variable of degree $-1$ and square zero. 
In this case, there is a generalization of  the Chern-Weil map constructed in [NT] :
 $$
{\operatorname{CW}}:{\mathbb H}^{\bullet}({\mathfrak h}_0[\epsilon], {\mathfrak h}_0; L^{\bullet}) \to {\mathbb H}^{\bullet}(\Sigma,{\mathcal L}^{\bullet}),
\leqno (4.4.3)
$$
which gives (4.4.1) when $L={\mathbb C}$ with the trivial action. Consider the following 
$({\widehat{\mathfrak g}}, H)$-modules:
$$
{\mathcal D}=\biggl\{\sum P_{\alpha}(x_1, \ldots, x_d, y_1, \ldots, y_k){\partial_x ^{\alpha}}\biggr\},
$$
where $P_{\alpha}$ are $r\times r$ matrices whose entries are power series, and
$$
\Omega ^{\bullet} =\biggl\{ \sum _I P_I(x_1, \ldots, x_d, y_1, \ldots, y_k) d^Ix\biggr\},
$$
which is  the space of differential forms whose coefficients are formal power series. The latter is a complex with the (fiberwise) De Rham differential. 
Moreover,  $\Omega ^{\bullet}$ is homotopy constant ($\epsilon {\widehat {\mathfrak g}}$ acts by exterior multiplication). 
The Hochschild, cyclic, etc. complexes of ${\mathcal D}$ inherit the $({\widehat{\mathfrak g}}, H)$-module structure; moreover, they also 
become homotopy constant (the $\epsilon X \in \epsilon {\widehat {\mathfrak g}}$
 acts by operators $\iota _X$ from  (4.2.4)). One constructs 
([BNT], pt. II, Lemma 3.2.4) 
a class
$$
\nu \in {\mathbb H}^0 ({\mathfrak h}_0[\epsilon], {\mathfrak h}_0; {\underline {\operatorname{Hom}}}
({\operatorname{CC}}^{\operatorname{per}}_{-{\bullet}}({\mathcal D}),\Omega^{2d+\bullet})), \leqno (4.4.3)
$$
such that ${\operatorname{CW}}(\nu)$ coincides with
 $$
\nu_{\mathcal D} \in {\mathbb H}^0 (X; {\underline {\operatorname{Hom}}}({\operatorname{CC}}^{\operatorname{per}}_{-{\bullet}}({\mathcal D}_{\Sigma/B}),
\Omega^{2d+\bullet}_{\Sigma/B})).
$$
To be precise, the cited lemma concerns the Weyl algebra of power series in both coordinates and derivations with the
Moyal product (clearly, differential operators of finite order form a subalgebra). Second, the construction there is for
the relative cohomology of the pair $(\Gg, {\mathfrak h})$ but it extends to the case of the pair $(\Gg[\epsilon], {\mathfrak h})$
of which $({\mathfrak h}_0 [\epsilon], {\mathfrak h}_0)$ is a sub-pair.

The cochain $\nu$ is actually independent of $y$. There is the canonical class
 $1$ in ${\operatorname {HC_0}}^{\operatorname {per}}({\mathcal D})$; it is ${\mathfrak h}_0$-invariant, and it is shown in [NT] how to extend it to a class in 
 ${\mathbb H}^0 ({\mathfrak h}_0[\epsilon], {\mathfrak h}_0; {\operatorname{CC}}^{\operatorname{per}}_{-{\bullet}}({\mathcal D}))$. On the other hand, 
$$
{\mathbb H}^0 ({\mathfrak h}_0[\epsilon], {\mathfrak h}_0; \Omega^{\bullet})
$$
can be naturally identified with
$$
{\mathbb H}^0 ({\mathfrak h}_0[\epsilon], {\mathfrak h}_0; {\mathbb C})
$$
It remains to show that 
$$
\nu(1) = \sum [{\operatorname{ch}}\cdot{\operatorname{Td}}]_{2(d+i)}\cdot u^i
$$
where ${\operatorname{ch}}$ is the corresponding invariant power series in $H^{\bullet}({\mathfrak {gl}}_r [\epsilon], 
{\mathfrak {gl}}_r; {\mathbb C})$
 and ${\operatorname{Td}}$ is the corresponding invariant power series in $H^{\bullet}({\mathfrak {gl}}_d [\epsilon],
  {\mathfrak {gl}}_d; {\mathbb C})$. 
 This was carried out in [BNT], Lemma 5.3.2.

\vfill\eject

\centerline{\bf 5. Comparison with the gerbe picture  }

\vskip 1cm
 
\noindent {\bf (5.1) $L^2$-sections of a vector bundle on a circle.}
Let $\Sigma$ be an oriented $C^\infty$-manifold diffeomorphic to the circle $S^1$
with the standard orientation, and let $E$ be a complex $C^\infty$-vector bundle on
$\Sigma$. Choose a smooth Riemannian metric $g$ on $\Sigma$ and a smooth Hermitian
metric $h$ on $E$. Let $\Gamma(\Sigma, E)$ be the space of $C^\infty$-sections of $E$. The choise
of $g,h$ defines a positive definite scalar product on this space and we denote by
$L^2_{g,h}(\Sigma, E)$ the Hilbert space obtained by completion with respect to this scalar product.

\proclaim (5.1.1) Lemma. For a different choice $g', h'$ of metrics on $\Sigma, E$
we have a canonical identification of topological vector spaces
$$L^2_{g,h}(\Sigma, E)\to L^2_{g', h'}(\Sigma, E).$$

\noindent{\sl Proof:} The Hilbert norms on $\Gamma(\Sigma, E)$ associated to $(g,h)$ and
$(g', h')$ are equivalent, since $\Sigma$ is compact.

\vskip .1cm

So we will denote the completion simply by $L^2(\Sigma, E)$. 

\vskip .2cm

Consider now the case when $\Sigma=S^1$ is the standard circle and $E=\bbC^r$ is the trivial bundle of rank $r$.
In this case $L^2(\Sigma, E) = L^2(S^1)^{\oplus r}$. Let us denote this Hilbert space by $H$. 
It comes with a polarization in the sense of Pressley and Segal [PS]. In other words, $H$ is decomposed
as $H_+\oplus H_-$ where $H_+, H_-$ are infinite-dimensional ortogonal closed subspaces
defined as follows. 

$H_+$ consists of vector-functions extending holomorphically into the unit disk $D_+ = \{ |z| <1\}$.
The space $H_-$ consists of vector functions extending holomorphically into the opposite annulus
$D_- = \{|z|>1\}$ and vanishing at $\infty$. 

The decomposition $H= H_+\oplus H_-$ yields the groups $GL_{res}(H)\subset GL(H)$, see
[PS] (6.2.1), as well as the Sato Grassmannian $Gr(H)$ on which $GL_{res}(H)$ acts transitively.
We recall that $Gr(H)$ consists of closed subspaces $W\subset H$ whose projection to $H_+$ is a 
Fredholm operator and the projection to $H_-$ is a Hilbert-Schmidt operator, see [PS] (7.1.1).

Given arbitrary $\Sigma, E$ as before, we can choose an orientation preserving diffeomorphism $\phi: 
S^1\to\Sigma$ and a trivialization $\psi: \phi^*E\to \bbC^r$. This gives an identification
$$u_{\phi, \psi}: L^2(\Sigma, E)\to H= L^2(S^1)^{\oplus r}.$$
In particular, we get a distinguished set of subspaces in $L^2(\Sigma, E)$, namely
$$Gr_{\phi, \psi}(\Sigma, E) = u_{\phi, \psi}^{-1}(Gr(H)),$$
and a distinguished subgroup of its automorphisms, namely
$$GL_{res}^{\phi, \psi}( L^2(\Sigma, E)) = u_{\phi, \psi}^{-1} GL_{res}(H) u_{\phi, \psi}.$$

\proclaim (5.1.2) Lemma. The subgroup $GL_{res}^{\phi, \psi}(L^2(\Sigma, E))$ and the
set $Gr_{\phi, \psi}(L^2(\Sigma, E))$ are independent on the choice of $\phi$ and $\psi$.

\noindent {\sl Proof:} Any two choices of $\phi, \psi$ differ by an element of the Atiyah group 
$AT(S^1, \bbC^r)$, see Example 3.3.5(c). This group being a semidirect product of
$\operatorname{Diffeo}(S^1)$ and $GL_r C^\infty (S^1)$, our statement follows from the known
fact that both of these groups are subgroups of $GL_{res}(H)$, see [PS]. 

So we will drop $\phi, \psi$ from the notation, writing $Gr(L^2(\Sigma, E))$ and
$GL_{res}(L^2(\Sigma, E))$. 

Recall further that $Gr(H)\times Gr(H)$ is equipped with a line bundle $\Delta$
(the relative determinantal bundle) which has the following additional structures:

\vskip .2cm

\noindent (a) Equivariance with respect to $GL_{res}(H)$.

\vskip .1cm

\noindent (b) A multiplicative structure, i.e., an identification
$$p_{12}^*\Delta\otimes p_{23}^*\Delta \to p_{13}^*\Delta\leqno (5.1.3) $$
of vector bundles on $Gr(H)\times Gr(H)\times Gr(H)$, which is
equivariant under $GL_{res}(H)$ and satisfies the associatiivity, unit and
inversion properties. 

\vskip .2cm

It follows from the above, that we have a canonically defined line bundle
(still denoted $\Delta$) on $Gr(L^2(\Sigma, E))$ equivariant under $GL_{res}(L^2(\Sigma, E))$
and equipped with a multiplicative structure. For $W, W'\in Gr(L^2(\Sigma, E))$ we denote
by $\Delta_{W, W'}$ the fiber of $\Delta$ at $(W, W')$.

As well known, the multiplicative bundle $\Delta$ gives rise to a category
($\bbC^*$-gerbe) $\cD et\, L^2(\Sigma, E)$ whose objects for the set
$Gr(L^2(\Sigma, E))$, while
$$\Hom_{\cD et\, L^2(\Sigma, E)} (W, W') = \Delta_{W, W'}-\{0\}.$$
The composition of morphisms comes from the identification
$$\Delta_{W, W'}\otimes \Delta_{W', W''}\to \Delta_{W, W''}$$
given by (5.1.3). 

\vskip .3cm

\noindent {\bf (5.2) $L^2$-direct image in a circle fibration.} Let now
$q:\Sigma\to B$ be a fibration in oriented circles and $E$ be a vector bundle
on $\Sigma$. We have then a bundle of Hilbert spaces $q_*^{L^2}(E)$ whose
fiber at $b\in B$ is $L^2(\Sigma_b, E_b)$. Further, by (5.2) this bundle has a
$GL_{res}(H)$-structure, where $H= L^2(S^1)^{\oplus r}$. Therefore we have the
associated bundle of Sato Grassmannians $Gr(q_*^{L^2}(E))$ on $B$ and the
(fiberwise) multiplicative line bundle $\Delta$ on
$$Gr(q_*^{L^2}(E))\times_B Gr(q_*^{L^2}(E)).$$
We define a sheaf of $C_B^{\infty*}$-groupoids on $B$ whose local objects
are local sections of $Gr(q_*^{L^2}(E))$ and for any two such sections defined
on $U\subset B$
$$\shHom (s_1, s_2) = (s_1, s_2)^*\Delta - 0_U,$$
where $0_U$ stands for the zero section of the induced line bundle.  This sheaf of groupoids
is locally connected and so gives rise to a $C_B^{\infty *}$-gerbe which we denote
$\cD et (q_* E)$. So we have the class
$$\bigl[ \cD et(q_* e) \bigr] \in H^2(B, C_B^{\infty *}).$$
Alternatively, consider the Atiyah group $G= AT(S^1, \bbC^r)$, see Example 3.3.5(c). 
By the above, $G\subset GL_{res}(H)$. The determinantal $\bbC^*$-gerbe $\cD et(H)$
(over a point) with $G$-action
 gives  a central extension $\widetilde{G}$ of $G$ by $\bbC^*$. A circle
fibration $q:\Sigma\to B$ gives a principal $G$-bundlle $P(\Sigma/B)$, as in (3.4.5),
and the following is clear.

\proclaim (5.2.1) Proposition. The gerbe $\cD et(q_*E)$ is equivalent to
 $\Lift_G^{\widetilde{G}}(P(\Sigma/B, E))$,
see Example 1.10.3. 

\vskip .2cm

Consider the exponential sequence (1.10.1)  of sheaves on $B$ and 
the corresponding coboundary map $\delta_2$, see (1.10.2).  
Then we have the class
$$\delta_2 \bigl[ \cD et(q_*E)\bigr] \in H^3(B, \bbZ).$$

\proclaim (5.2.2) Theorem. The image of $\delta \bigl[ \cD et(q_*E)\bigr]$
in $H^3(B, \bbC)$ coincides with negative of the class $C_1(q_*E)$ defined
in (3.5.2).  

To prove Theorem 5.2.2, we apply Proposition 3.7.7 
  to $G= AT(S^1, \bbC^r)$ and $\beta$ being the class of the  central extension $\widetilde{G}$.
  Then $\Gg = \cA_{\bbC^r}(S^1)$ is the
Atiyah algebra of the trivial bundle on $S^1$ and $\gamma$ is the class of the ``trace" central extension
induced from the Lie algebra ${\frak{gl}}_{res}(H)$ of $GL_{res}(H)$. We have the embeddings
$$\Gg \subset {\frak{gl}}_r(\cD(S^1)) \subset {\frak{gl}}_{res}(H),$$
and the trace central extension is represented by an explicit cocycle $\Psi$  of ${\frak{gl}}_{res}(H)$ 
(going back to [T]). 
Let $z$ be the standard complex coordinate on $S^1$ such that $|z|=1$. Then the formula for the restriction of
$\Psi$ to ${\frak{gl}}_r(\cD(S^1))$ was given in [KP], see also [KR], formula (1.5.2): 
$$\Psi(f(z) \partial_z^m, g(z) \partial_z^n) = {m! n!\over (m+n+1)!} \Res_{z=0} dz \cdot \Tr (f^{(n+1}(z) g^{(m)}(z)), 
\leqno (5.2.3)$$
where $f^{(n)}$ means the $n$th derivative in $z$. 
Our statement now reduces to the following.

\proclaim (5.2.4) Lemma. The second Lie cohomology class of ${\frak{gl}}_r \cD(S^1)$ given by the cocycle $\Psi$
is equal to the negative of the class corresponding to the fundamental class of $S^1$ via the identification (2.3.4).

\noindent {\sl Proof:} As the space of (continuous) Lie algebra homology in question is 1-dimensional,
it is enough to evaluate the cocycle $\Psi$ on the Lie algebra 2-homology class $\sigma$ from (2.3.4) and to show that this
value is precisely equal to 1. For this it is enough to consider $r=1$. Let $\cD = \cD(S^1)$ 
for simplicity. 

We need to recall the explicit form  of the identification (2.3.1) for the case $n=1$ (first Hochschild
homology maps to the second Lie algebra homology). In other words, we need to recall the definition of
the map.
$$\epsilon: HH_1(\cD)\to H_2^{\Lie}({\frak{gl}}(\cD))\to \bbC. $$
As explained in [BG] and [W], this map is defined via the order filtration $F$ 
on the ring $\cD$ and uses the corresponding spectral sequence.  
This means we need to start with a Hochschild
1-cycle $\sigma = \sum P_i\otimes Q_i\in \cD\otimes\cD$ and form its highest symbol cycle
$$\Smbl(\sigma) = \sum \Smbl(P_i)\otimes \Smbl(Q_i)\in \gr(\cD)\otimes\gr(\cD),$$
which gives an element in $\Hoch_1(\gr(\cD))$. As $\gr(\cD)$ is the ring of
polynomial functions on $T^*S^1$, Hochschild-Kostant-Rosenberg gives
$HH_1(\gr(\cD)) = \Omega^1({T^*S^1})$, the space of 1-forms on $T^*S^1$ polynomial
along the fibers. So the class of $\Smbl(\sigma)$ is a 1-form $\omega = \omega(\sigma)$
in $T^*S^1$. This is an element of the $E_1$-term of the spectral sequence for the Hochschild
homology of the filtered ring $\cD$.

Further, one denotes by $*$ the symplectic Hodge operator in forms on $T^*S^1$. 
The results  of {\em loc. cit.} imply the differential in the $E_1$-term is $*d*$ where $d$ is the
de Rham differential on $T^*S^1$ while higher differentials vanish. This means that under our assumptions
 $*\omega(\sigma)$ is a closed 1-form and
$$\epsilon(\sigma) = \int_{S^1} *\omega(\sigma).$$
To show Lemma 5.3.5 we need to exhibit just one $\sigma$ as above such that
 $$0\neq \epsilon (\sigma) = \Psi(\sigma) := \sum \Psi(P_i, Q_i) . $$

We take 
$$\sigma = z^2\otimes z^{-1} \partial_z - 2 z\otimes\partial_z.$$
Then one sees that $\sigma$ is a Hochschild 1-cycle and  $\Psi(\sigma) = 1$. 
On the other hand, let $\theta$ be the real coordinate on $S^1$ so that
$z=\exp(2\pi i \theta)$. Then the real coordinates on $T^*S^1$ are $\theta, \xi$
with $\xi = \Smbl(\partial/\partial \theta)$, so the Poisson bracket
$\{\theta, \xi\}$ is equal to 1. In terms of the coordinate $z$ it means
that $\xi = \Smbl(z \partial/\partial z)$ and $\{z, \xi\} = z$. 
Therefore
$$\Smbl (\sigma) = z^2\otimes z^{-2}\xi - 2 z\otimes z^{-1}\xi$$
and hence
$$\omega (\sigma) = z^2 d(z^{-2}\xi) = 2 z d(z^{-1}\xi) = -dz - z^{-1}\xi,$$
see [L] p. 11.
The symplectic (volume) form on $T^*S^1$ is $(dz/z)\wedge d\xi$, so the symplectic Hodge
operator is given by
$$*d\xi = dz/z, \quad *dz/z = d\xi, \quad *^2=1.$$
Therefore
$$*\omega(\sigma) = -dz/z -\xi d\xi, \quad \int_{S^1} *\omega(\sigma) = -1$$
 and we are done.

\vfill\eject

\centerline{\bf References}

\vskip 1cm

\noindent [AH] M.~F.~ Atiyah, F.~ Hirzebruch, {\em Riemann-Roch theorems for
differentiable manifolds}, Bull. Amer. Math. Soc. {\bf 65} (1959), p. 276--281. 

\vskip .1cm

\noindent [Bl] S.~Bloch, $K_2$ {\em and algebraic cycles},  Ann. of Math. (2)  {\bf 99}  (1974), 349--379. 

\vskip .1cm

\noindent [Bo] R. Bott, {\em On the characteristic classes of groups of diffeomorphisms},
L'Enseignment Math. {\bf 23} (1977), 209-220.

\vskip .1cm

\noindent [BNT] P.~Bressler, R.~Nest, B.~Tsygan, {\em Riemann-Roch theorems via deformation quantization. I, II}, 
 Adv. Math. {\bf 167}  (2002),  no. 1, p. 1--25, 26--73. 
 
 \vskip .1cm
 
 \noindent [Bry] L.~ Breen, {\em On the classification of 2-gerbes and 2-stacks,} Ast\'erisque,
   {\bf  225}  (1994), 160 pp.
 
 \vskip .1cm
 
 \noindent [Bry] J.~ L.~ Brylinski, {\em Loop Spaces, Characteristic Classes and Geometric
 Quantization}, Birkhauser, Boston, 1993.

 \vskip .1cm
 
 \noindent [BG] J.-L.~ Brylinski, E. Getzler, {\em  The homology of algebras of pseudodifferential symbols 
 and the noncommutative residue}, 
$K$-Theory {\bf 1} (1987), no. 4, 385--403.

 \vskip .1cm

 \noindent [C] K.~T.~Chen, {\em Iterated path integrals},  Bull. Amer. Math. Soc.  {\bf 83}  (1977), no. 5, 831--879. 
 
 \vskip .1cm
 
 \noindent [D1] P. Deligne, {\em Le determinant de la cohomologie}, 
  Current trends in arithmetical algebraic geometry (Arcata, Calif., 1985),  93--177, Contemp. Math., 
  67, Amer. Math. Soc., Providence, RI, 1987.
 
 \vskip .1cm
 
 \noindent [D2] P.~ Deligne, {\em Le groupe fondamental de la droite projective moins trois points}, 
  Galois groups over $Q$ (Berkeley, CA, 1987),  79--297, Math. Sci. Res. Inst. Publ., {\bf 16}, Springer, New York, 1989. 
 
 \vskip .1cm
 
 \noindent [E] R.~ Elkik, {\em Fibr\'es d'intersections et int\'egrales de classes de Chern,}   
   Ann. Sci. \'Ecole Norm. Sup. (4)  {\bf 22}  (1989),  no. 2, 195--226. 

\vskip .1cm

\noindent [FT] B.~Feigin, B.~Tsygan, {\em Riemann-Roch theorem and Lie algebra cohomology, I},
 Proceedings of the Winter School on Geometry and Physics (Srn\'{i}, 1988). 
 Rend. Circ. Mat. Palermo (2) Suppl.  {\bf  21} (1989), 15--52.

\vskip .1cm

\noindent [F] D.~ B.~ Fuks, {\em Cohomology of Infinite-Dimensional Lie Algebras}, Consultants Bureau, New York and London, 1986.

\vskip .1cm

\noindent [HS] V.~Hinich, V.~Schechtman, {\em Deformation theory and Lie algebra homology, I,II},
 Algebra Colloq. {\bf 4} (1997), no. 2, 213--240, 291--316.
 
 \vskip .1cm
 
 \noindent [KP] V.~ Kac, D.~ Peterson, {\em Spin and wedge representations of infinite-dimensional Lie algebras
 and groups}, Proc. Natl. Acad. Sci. USA, {\bf 78} (1981), 3308-3312. 
 
 \vskip .1cm
 
 \noindent [KR] V.~ Kac, A.~Radul, 
 {\em Quasifinite highest weight modules over the Lie algebra of differential operators on the circle,}
   Comm. Math. Phys.  {\bf 157}  (1993),  no. 3, 429--457. 
 
 \vskip .1cm
 
 \noindent [Kal] R.~ Kallstrom, {\em Smooth Modules over Lie Algebroids I}, preprint math.AG/9808108.

\vskip .1cm

\noindent [KV1] M.~Kapranov, E.~Vasserot, {\em Vertex algebras and the formal loop space},
 Publ. Math. Inst. Hautes Etudes Sci. {\bf 100} (2004), 209--269. 
 
 \vskip .1cm
 
 \noindent [KV2] M.~Kapranov, E.~Vasserot, {\em Formal loops II:  a local Riemann-Roch
 theorem  for determinantal gerbes}, preprint math.AG/0509646.  
 
 \vskip .1cm
 
 \noindent [LM] G.~Laumon, L.~ Moret-Bailly, {\em Champs Alg\'ebriques}, 
   Springer-Verlag, Berlin, 2000.

\vskip .1cm

\noindent [L] J.~L. Loday, {\em Cyclic Homology}, Second edition,
   Springer-Verlag, Berlin, 1998. 

\vskip .1cm

\noindent [Lo]  J.~Lott, {\em Higher-degree analogs of the determinant line bundle},
  Comm. Math. Phys. {\bf 230}  (2002),  no. 1, 41--69.
  \vskip .1cm
  
  \noindent [Mac] K.~C.~H.~McKenzie, {\em General Theory of Lie Groupoids and Lie Algebroids},
   London Mathematical Society Lecture Note Series, {\bf 213}. Cambridge University Press, Cambridge, 2005.
   
   \vskip .1cm
   
   \noindent [Mo] S.~ Morita, {\em  Geometry of Characteristic Classes},   
    Translations of Mathematical Monographs, {\bf 199}. Iwanami Series in Modern Mathematics.
     American Mathematical Society, Providence, RI, 2001
  
  \vskip .1cm

\noindent [NT] R.~Nest, B.~Tsygan, {\em Algebraic index theorem for families},
 Adv. Math. {\bf  113} (1995),  no. 2, 151--205. 
 
 \vskip .1cm
 
 \noindent [PS] A. Pressley, G.~B.~ Segal, {\em Loop Groups}, Cambridge University Press, 1986.
 
 \vskip .1cm
 
 \noindent [RSF] A.~G.~ Reiman, M.~ A.~ Semenov-Tyan-Shanskii, L.~ D.~ Faddeev, 
 {\em Quantum anomalies and cocycles on gauge groups,} Funkt. Anal. Appl. {\bf 18} (1984), No. 4, 64-72. 
 
 \vskip .1cm
 
 \noindent [T] J.~ Tate, {\em Residues of differentials on curves}, 
  Ann. Sci. École Norm. Sup. (4) {\bf  1 } (1968) 149--159. 

\vskip .1cm

\noindent [W]  M.~Wodzicki, {\em Cyclic homology of differential operators},
 Duke Math. J.  {\bf 54}  (1987),  no. 2, 641--647. 

\vskip 1cm

\end{document}